\documentclass[12pt]{article}

\scrollmode
\usepackage{color}
\usepackage{amsfonts}
\usepackage{amsmath}
\usepackage{amssymb}
\usepackage{graphicx}

\usepackage{mathptmx}

\def\eps{\varepsilon}
\def\cinf{\underline c}
\def\qed{\hfill\rule{.2cm}{.2cm}}
\def\P{{\mathbb P}}
\def\E{{\mathbb E}}
\def\Z{{\mathbb Z}}
\def\N{{\mathbb N}}
\def\R{{\mathbb R}}
\def\Q{{\mathbb Q}}

\def\1{{\mathbf 1}}

\def\a{\alpha}
\def\o{\omega}  
\def\l{\lambda} \def\hl{\hat \l}
\def\g{\gamma}

\def\bB{\bar B}

\newcommand{\AAA}          {\mathcal{A}}

\newcommand{\BB}          {\mathcal{B}}
\newcommand{\EE}          {\mathcal{E}}
\newcommand{\FF}          {\mathcal{F}}

\newcommand{\LL}         {\mathcal{L}}
\newcommand{\MM}         {\mathcal{M}}
\newcommand{\HH}         {\mathcal{H}}
\newcommand{\WW}         {\mathcal{W}}

\def\O{\Omega}

\newtheorem{theorem}{Theorem}[section]
\newtheorem{lemma}[theorem]{Lemma}

\newtheorem{proposition}[theorem]{Proposition}
\newtheorem{definition}[theorem]{Definition}
\newtheorem{remark}[theorem]{Remark}

\def\beq{\begin{equation}}
\def\eeq{\end{equation}}
\newcommand{\bei}{\begin{itemize}}
\newcommand{\eei}{\end{itemize}}
\newcommand{\ben}{\begin{enumerate}}
\newcommand{\een}{\end{enumerate}}
\newcommand{\beqn}{\begin{eqnarray}}
\newcommand{\beqnn}{\begin{eqnarray*}}
\newcommand{\eeqn}{\end{eqnarray}}
\newcommand{\eeqnn}{\end{eqnarray*}}
\newcommand{\brm}{\begin{rmk}}
\newcommand{\erm}{\end{rmk}}

\newcommand {\tr}[1] {{#1}\cdot }

\begin{document}

\title{Einstein relation for reversible diffusions in random environment }
\author{
N.~Gantert  \footnote{Technische Universit\"at M\"unchen,
Fachbe\-reich Mathematik und Informatik,
Boltzmannstr. 3, 85748 Garching,
GERMANY. gantert@ma.tum.de
}
\and
P.~Mathieu \setcounter{footnote}{6}\footnote{ Universit\'e de Provence, CMI, 39 rue Joliot-Curie, 13013 Marseille, FRANCE.
pierre.mathieu@cmi.univ-mrs.fr}
\and
A.~Piatnitski \setcounter{footnote}{3}\footnote{Lebedev Physical Institute of Russian Academy of Sciences and
Narvik Institute of Technology, P. O. Box 385, N-8505 Narvik, NORWAY.
andrey@sci.lebedev.ru}
}

\maketitle

\begin{abstract}
We consider reversible diffusions in random environment and prove the Einstein relation for this model.
It says that the derivative at $0$ of the effective velocity  under an additional local drift equals
the diffusivity of the model without drift. The Einstein relation is conjectured to hold for a variety of models but it is proved insofar only in
particular cases. Our proof makes use of homogenization arguments, the Girsanov transform, and a refinement of the regeneration times introduced in \cite{kn:LS}.
\end{abstract}




\section{Introduction}

The present paper deals with diffusions in a random stationary environment, a model
for the dynamics of particles in a disordered medium at thermal
equilibrium. This subject has been the object of intense research over the
past thirty years. In spite of many recent progresses, see \cite{kn:SzZeitouni} for
instance, many questions regarding the long time behavior of these
processes, such as laws of large numbers or central limit theorems, are
still challenging open problems.

We shall only be concerned with reversible dynamics. In this context the
idea of the environment seen from the particle, as discussed in \cite{kn:KV} or
\cite{kn:DFGW},
provides a powerful tool to adapt the 'corrector approach' from
homogenization theory and eventually prove invariance principles. One then
shows that the trajectory of a particle evolving in such an environment,
in a large time scale, behaves like a Brownian motion with mean square
displacement proportional to time, the proportionality being expressed by
the asymptotic covariance or {\it effective diffusivity} matrix $\Sigma$. A good
understanding of $\Sigma$ is thus of primary interest.

>From reversibility follows a variational formula for the effective diffusivity,
see \cite{kn:Ko85} and \cite{kn:JKO} for the discrete and continuous cases, respectively.
Many works in PDE or theoretical physics address the question of estimating the effective diffusivity.
See \cite{kn:FM} or \cite{kn:FSS} for instance.
Here we provide a completely different interpretation of $\Sigma$ as the so-called {\it mobility}.

In the series of papers A. Einstein devoted to Brownian motion at
the beginning of past century, see \cite{kn:Einstein}, along with
the diffusivity matrix, the great physicist introduced another
important quantity called the {\it mobility}. The mobility measures
the response of the diffusing particle to a constant exterior force.
Let us consider the perturbed process obtained by imposing a
constant drift of strength $\lambda$ in some fixed direction. One
would expect the perturbed process to satisfy a law of large numbers
with effective drift  $\ell(\lambda)$. The mobility can then be
interpreted as the derivative of $\ell(\lambda)$ as $\lambda$ tends
to $0$. Einstein claimed that the mobility and the diffusivity of a
particle coincide.

This 'diffusivity vs mobility' relation played a central role in
Einstein's theory of molecular diffusivity because it was amenable
to experimental verifications and eventually lead to evidences that
matter is molecular. Since then, Einstein's relation opened the way
to important developments both in experimental physics, with J.
Perrin's Nobel prize, and theoretical physics, with connections to
the Fluctuation Dissipation Theorem or the Green-Kubo relation.
Authors in these fields usually assume the Einstein relation, on the
basis on heuristics, without being actually able to prove it, see
\cite{kn:FSS} for instance. Indeed, only very few rigorous papers
investigate the mathematical contents of the Einstein relation and
finding a general strategy for proving it remained an open
mathematical problem for years.
 Note that the mere existence of the mobility i.e. the fact that
$\ell(\lambda)$ has a derivative at $\lambda=0$ is far from obvious.

In the case of  periodic environments the Einstein relation can be easily checked by direct computations which are mostly reduced to differentiating the perturbed cell problem on the torus. More generally, if the process of the
environment seen from the particle has a spectral gap, the Einstein relation can be proved by perturbation theory, see
 \cite{kn:KO2}.

The first consistent mathematical approach to the Einstein relation
for random environments was attempted by J. Lebowitz and H. Rost in
\cite{kn:LR}. These authors then avoided the difficulty of proving
the existence of the mobility by considering exterior forces whose
strength vanishes as time goes to infinity in such a way that the
perturbed process still has a limit on the diffusive scale. The
mobility is then defined as the asymptotic mean displacement of the
particle; no derivative is involved and this weak form of the
Einstein relation does not tell us anything about the effective
drift induced by a constant exterior force. On the other hand, in
this weaker form the Einstein relation holds in the general ergodic
reversible case. We shall discuss the Lebowitz-Rost theory in more
details later, see the end of section \ref{II.4}.

More recently, in \cite{kn:KO}, T. Komorowski and S. Olla proved  the
Einstein relation for random walks with random conductances in dimension
higher than $3$ and with the extra restriction that conductances are only
allowed to take two different positive values. Their approach, which is an adaptation of \cite{kn:loulakis} 
to the RWRE case,
 is quite
different from ours  and it is not clear
how it can be adapted to more general models, as for instance
random walks on percolation clusters, see \cite{kn:BGP}
and \cite{kn:Sz}.

In the present paper, we shall prove the Einstein relation for symmetric
diffusions driven by  a random environment with bounded potential and short range correlations.
A more detailed definition of the model and precise assumptions are given
in part \ref{sec:II}. As discussed above, the invariance principle, and therefore the
existence of the diffusivity, were already known for such models in the
eighties. The law of large numbers for the perturbed dynamics, i.e. the
existence of $\ell(\lambda)$, is ensured by a more general theorem of L. Shen
in \cite{kn:LS}
on renewal properties of drifted diffusions in random environments. The
key issue in proving the Einstein relation is then to show that $\ell(\lambda)$
has a derivative at $\lambda=0$ and to identify this derivative with the effective diffusivity
of the non-perturbed diffusion. To this end we
compute the mean
position of the perturbed diffusion on the critical scale $\lambda^2 t=1$.
This is indeed the scale on which the diffusive and drift components of
the perturbed dynamics equilibrate. Thus we obtain two
different expressions for the mean displacement in terms of either the
diffusivity matrix $\Sigma$, see equation (\ref{eq:***}), or in terms of $\ell(\lambda)$, see (\ref{eq:lln}).
Identifying
both expressions gives the Einstein relation.

The proofs combine different ingredients: homogenization arguments
and Girsanov transforms, - see Section \ref{sec:III} - PDE estimates
and a-priori bounds on hitting times for perturbed diffusions - see
Section \ref{sec:IV} - and renewal arguments - see Section
\ref{sec:V}. All these ingredients had already appeared in the
literature but, in order to treat the critical scale $\lambda^2
t=1$, we had to refine many arguments and often introduce
alternative  strategies as, for instance, with the regeneration
times in Section~\ref{sec:V}.

In Sections \ref{sec:II} - \ref{sec:VI}, we focus only on smooth
environments which allows us to use stochastic differential
equations. In Section \ref{sec:VII}, we relax this smoothness
assumption and treat the case of (still bounded) but only measurable
environments, relying on Dirichlet form theory.

One might hope that our approach could be adapted to apply to other
models of diffusions or random walks in random environments.


\vskip 1cm
\section{Model and statement of the theorem}
\setcounter{equation}{0}
\label{sec:II}

\subsection{Diffusions in a random environment} \label{ssec:II.1}

We shall be dealing with diffusion processes in $\R^d$ whose generators are
of the form
\beqn \label{eq:gen} \LL^\o f(x)=\frac 1 2 e^{2V^\o(x)}\, \mathrm{div}(e^{-2V^\o}\,a^\o\,\nabla f)(x)\,,
\eeqn
where $a^\o$ and $V^\o$ are realizations of a random environment with finite
range of dependence.

More precisely, our assumptions are as follows.

Let $(\Omega,\AAA,\Q)$ be a probability space equipped with a group
action of $\R^d$ 
that we denote with $(x,\o)\rightarrow x.\o$.
We also assume that the map $(x,\o)\rightarrow x.\o$ is $(\BB_d\times\AAA,\AAA)$-measurable,
where $\BB_d$ is the Borel $\sigma$-field on $\R^d$.

\medskip\noindent
{\bf Assumption 1:}  the action $(x,\o)\rightarrow x.\o$ preserves the measure $\Q$ and is ergodic.

Let $V$ be a measurable real-valued function on $\O$ and let $\sigma$ be a measurable
function on $\Omega$ taking its values in the set of real $d\times d$ symmetric matrices. Define
$$ V^\o(x)=V(x.\o)  \,,\, \sigma^\o(x)=\sigma(x.\o)  \,.$$
We also introduce the notation
$$a^\o=(\sigma^\o)^2 \hbox{\, and\, } b^\o=\frac 12 \mathrm{div}a^\o-a^\o\nabla V^\o\,.$$

Observe that both $\sigma^\o$ and $b^\o$ are then stationary fields i.e.
$\sigma^\o(x)=\sigma(x.\o)$ and $b^\o(x)=b(x.\o)$ for some functions $\sigma$ and $b$.

\medskip\noindent
{\bf Assumption 2:} for any environment $\o$, the functions $x\rightarrow V^\o(x)$ and
$x\rightarrow \sigma^\o(x)$ are smooth.
To avoid triviality, we also assume that at least one of them is not constant.

\medskip\noindent
{\bf Assumption 3:} $V$ is bounded and  $a^\o$ is uniformly elliptic, namely there exists a constant $\kappa$
such that, for all $\o$, $x$ and $y$,
\beqn\label{eq:ell} \kappa\vert y\vert^2\leq \vert \sigma^\o(x)y\vert^2\leq \kappa^{-1} \vert y\vert^2\,.
\eeqn

For a Borel subset $F\subset\R^d$, we define the $\sigma$-field
$$\HH_F=\sigma\{V(x.\o),\sigma(x.\o)\,:\, x\in F\}$$
and we assume the following independence condition:

\medskip\noindent
{\bf Assumption 4:} there exists $R$ such that for any Borel subsets $F$ and $G$
such that $d(F,G)>R$ (where $d(F,G)=\inf\{\vert x-y\vert\,:\, x\in F, y\in G\}$ is the distance between $F$ and $G$) then
\beqn \label{eq:indep} \HH_F \hbox{ \rm and } \HH_{G} \hbox{ are independent}\,.
\eeqn

Let $(W_t:t\geq 0)$ be a Brownian motion defined on some probability
space $(\WW,\FF,P)$. We denote expectation with respect to $P$ by
$E$. By {\it diffusion in the environment $\o$} we mean the solution
of the stochastic differential equation \beqn \label{eq:sde1}
dX_x^\o(t)=b^\o(X_x^\o(t))\,dt+\sigma^\o(X_x^\o(t))\,dW_t \,;\,
X_x^\o(0)=x\,. \eeqn Then $X_x^\o$ is indeed the Markov process
generated by the operator $\LL^\o$ in equation (\ref{eq:gen}). We
shall denote with $P^\o_x$ the law of $X^\o$ on the path space
$C(\R_+,\R^d)$. It is usually referred to as the {\it quenched} law
of the diffusion in a random environment. We will also need the
so-called {\it annealed} law: \beqn\label{eq:ann} \P_x[A]: =\int
d\Q(\o)\int dP^\o_x(w)\1_A(\o,w)\,, \eeqn for any measurable subset
$A\subset \O\times C(\R_+,\R^d)$.

Expectation with respect to $P^\o_x$ will be denoted with $E^\o_x$ and expectation with respect to
$\P_x$ will be denoted with $\E_x$.

We use the notation $X(t)$ for the coordinate process on path space $C(\R_+,\R^d)$.

\vskip1cm
\subsection{ Effective diffusivity} \label{ssec:II.2}

\begin{definition}\label{df:Sigma}
Let $\Sigma$ be the {\it effective diffusivity} matrix defined by
\beqn\label{eq:Sigma}
e\cdot\Sigma\, e: =\lim_{t\rightarrow +\infty} \frac 1 t \E_0[(e\cdot X(t))^2]\,, \eeqn
where $e$ is any vector in $\R^d$ and $x\cdot y$ denotes the scalar product of the two
vectors $x$ and $y$.
\end{definition}

The fact that the limit in (\ref{eq:Sigma}) exists is (almost) a consequence of the Central Limit Theorem
for the process $X$ under $\P_0$. More is actually known:
$X$ satisfies a full invariance principle.
Namely: for almost any realization of the environment $\o$, the laws of the sequence of rescaled processes
$(X^\eps(t)=\eps X(t/\eps^2)\,;\, t\geq 0)$ under $P^\o_0$ weakly converge as $\eps$ goes to $0$ to the law of
a Brownian motion with covariance matrix $\Sigma$. References on this Theorem include
\cite{kn:DFGW}, \cite{kn:KV}, \cite{kn:Ko80}, \cite{kn:Osada}, \cite{kn:PV}
among others.
The convergence of the variance of the process to $\Sigma$ is explicitly stated in \cite{kn:DFGW}
formula (2.44).

The invariance principle also has a PDE counterpart in terms of homogenization theory,
see for instance the book \cite{kn:JKO}.
The generator of the process $X^\eps$ under $P^\o_x$ is the rescaled elliptic operator
with rapidly oscillating coefficients
$$\frac 12 a^\o(\frac .\eps)\Delta +\frac 1{2\eps}b^\o(\frac . \eps)\nabla\,.$$
Its limit, in the sense of homogenization theory, is the elliptic operator with constant
coefficient $$\frac 1 2 \mathrm{div}(\Sigma\nabla)\,,$$ where $\Sigma$
is the same matrix as in (\ref{eq:Sigma}).

The effective diffusivity $\Sigma$ is a symmetric matrix.
As a consequence of Assumption 1 on ergodicity, $\Sigma$ is deterministic (i.e. $\Sigma$ does not
depend on $\o$). Furthermore, due to the ellipticity Assumption 3, $\Sigma$
is also known to be positive definite.


In general, there is no simple expression of $\Sigma$ in terms of
$\sigma$ or $V$. (For instance, $\Sigma$ is by no means the average
of $\sigma^2$!)
The proof of the homogenization theorem actually provides an expression of $\Sigma$ as a function of
the solution of a Poisson equation - the so-called {\it corrector approach}. Since the operator
$\LL^\o$ is self-adjoint with respect to the measure $e^{-2V^\o(x)} dx$, there is also a
variational formulation of the Poisson equation and therefore a variational formula for $\Sigma$ itself.
We will not need it in this paper.
Our main Theorem actually gives a quite different interpretation of
$\Sigma$ as the {\it mobility} of $X^\o$, see below.

\vskip1cm
\subsection{Perturbed diffusions} \label{II.3}

We shall now consider perturbations of the process $X$ obtained by inserting
a local drift in equation (\ref{eq:sde1}).

We use the following notation. Let $e_1$ be a non-zero vector with $|e_1| = 1$ and $\lambda>0$. We define $\hat\lambda$
to be the vector ${\hat\lambda}=\lambda e_1$. We think of $e_1$ as being fixed while $\lambda$ is due to
tend to $0$.  We assume throughout the whole paper that $\lambda\leq 1$.

Let us consider the perturbed stochastic differential equation:
\beqn \label{eq:sde2}
dX_x^{\l,\,\o}(t)=b^\o(X_x^{\l,\,\o}(t))\,dt+
\,a^\o(X_x^{\l,\,\o}(t))\hl\,dt+\sigma^\o(X_x^{\l,\,\o}(t))\,dW_t \,;\, X_x^{\l,\,\o}(0)=x\,.
\eeqn

The process $X_x^{\l,\,\o}$ is now a Markov process with generator
\beqn\label{generator}
 \LL^{\l,\,\o} f(x)
&=&
\frac 1 2 e^{2V^\o(x)}\, \mathrm{div}(e^{-2V^\o}\,a^\o\,\nabla f)(x)+\,a^\o(x)\hl\,\cdot \nabla f(x)\cr
&=&\label{eq:gen'}\frac 1 2 e^{2V^{\l,\,\o}(x)}\, \mathrm{div}(e^{-2V^{\l,\,\o}}\,a^\o\,\nabla f)(x)\,,
\eeqn
where $V^{\l,\,\o}(x)=V^\o(x)-\tr{\hl}\, x$.
We shall use the notation $P_x^{\l,\,\o}$ for the law of $X_x^{\l,\,\o}$, $E_x^{\l,\,\o}$ for the corresponding
expectation as well as $\P_x^\l$ and $\E_x^\l$ for the annealed probability and expectation  defined analogously to
(\ref{eq:ann}).

Our model is a special case of diffusions with drifts considered by
L. Shen in \cite{kn:LS}  for which the author proved a law of large
numbers:
for almost any environment $\o$, the ratio $X(t)/t$ has an almost
sure limit under $P_0^{\l,\,\o}$, 
say $\ell(\l)$. The convergence also holds in $L^1(\P^\l_0)$. Moreover
$\ell(\l)$ is deterministic and $\tr{\hat\l}\ell(\l)>0$. Note that
the proof strongly relies on the independence property Assumption 4.
We thus define the effective velocity:

\begin{definition} \label{df:velo}
Let $\lambda>0$. Let $\ell(\l)$ be the {\it effective drift} vector defined by
\beqn\label{eq:vlambda}
\ell(\l)=\lim_{t\rightarrow +\infty} \frac 1 t \E_0^\l[X(t)]\,. \eeqn
By convention $\ell(0)=0$.
\end{definition}

\vskip1cm
\subsection{The Einstein relation} \label{II.4}

Our main result is the following theorem.

\begin{theorem}  \label{theo}
The function $\lambda\rightarrow \ell(\l)$
has a derivative at $\lambda=0$ which satisfies
\beqn\label{eq:einstein}
\lim_{\l\rightarrow 0} \frac 1 \l \,\ell(\l)=\Sigma\, e_1\,.\eeqn
\end{theorem}

This justifies the following definition.
\begin{definition}\label{mobdef}
Call {\it mobility} in the direction $e_2$
the derivative at $\l=0$ of the velocity $\tr{e_2}\, \ell(\l)$.
\end{definition}

Theorem \ref{theo} and Definition \ref{mobdef}  can be compared with the main result of \cite{kn:LR}
where the mobility is defined as
the mean position of the process $X^\eps$ under $\P_0^\eps$. The authors prove that, as $\eps$ tends to $0$,
the law of  $X^\eps$ under $\P_0^\eps$ converges to the law of Brownian motion with drift $v$
given by $v=\Sigma\,e_1$. These results are consequences of the invariance
principle under $\P_0$ and do not require any information on the asymptotic behaviour of the process
under $\P_0^\l$ for a fixed $\lambda$
(and indeed the law of large numbers of L. Shen was not known at the time \cite{kn:LR} was written).


\vskip 1cm
\section{Girsanov transforms}
\setcounter{equation}{0}
\label{sec:III}

The aim of this section is to establish Proposition \ref{lm:LR} below. In this part of the paper we only use Assumptions 1, 2, 3.

\begin{proposition} \label{lm:LR}
Let $\alpha \geq 1$. Then
\beqn \label{eq:***}
\lim_{\l\rightarrow 0\,;\, t\rightarrow +\infty \,;\, \l^2t=\alpha} \E_0^\l\left[\,\frac {X(t)}{\l t}\right]
= \Sigma\, e_1\,.\eeqn
Also
\beqn \label{eq:max****}
\sup_{\alpha\geq 1} \limsup_{\l\rightarrow 0\,;\, t\rightarrow +\infty \,;\, \l^2t=\alpha}
\E_0^\l\left[\,\frac {\max\limits_{s \le t}\vert X(s)\vert^2}{\l^2 t^2}\right]
<\infty\,,
 \eeqn

\end{proposition}

\begin{remark} Observe that (\ref{eq:max****}) directly follows from  Lemma \ref{lm:ext}
in the next section of the paper.
\end{remark}

\noindent
{\it Proof of Proposition \ref{lm:LR} in the case $V=0$}

We first prove Proposition \ref{lm:LR} in the case there is no potential i.e.  we start assuming that $V=0$.

We use Girsanov transforms pretty much as in \cite{kn:LR} i.e. the explicit expression of the Radon-Nikodym derivative of
$P_0^{\l,\,\o}$ with respect to $P_0^\o$.

Let us first recall Girsanov transforms, see  \cite{kn:RYor}, chapter VIII.
Let $X= (X_1, \ldots ,X_d)$ and $Y=(Y_1, \ldots ,Y_d)$  be solutions of stochastic differential
equations of the form
\beqnn
dX(t)=b(X(t))\,dt+\sigma(X(t))\,dW_t \,;\, X(0)=x\,,
\eeqnn
and
 \beqnn
dY(t)=c(Y(t))\,dt+\sigma(Y(t))\,dW_t \,;\, Y(0)=x\,,
\eeqnn
where $(W_t:t\geq 0)$ is a Brownian motion and the coefficients $\sigma$, $b$ and $c$ are subject
to smoothness and ellipticity assumptions as in Assumptions 2 and 3.
Let $P^X$ and $P^Y$ be the laws of the processes $X$ and $Y$ on the path space $C(\R_+,\R^d)$.
Let $\FF_t=\sigma\{ X(s)\,;\, s\leq t\}$ be the filtration generated by the coordinate process up to time $t$.
Then the restriction of $P^Y$ to $\FF_t$ is absolutely continuous with respect to the restriction of $P^X$
to $\FF_t$ and the Radon-Nikodym derivative is given by the Girsanov formula:
$$E[F(Y([0,t]))]=E[F(X([0,t]))\,e^{M(t)-\frac{1}{2}\langle M\rangle(t)}]\,,$$
for any time $t$,
for any bounded continuous function $F$ on $C([0,t],\R^d)$ and
where $M$ is the martingale
\beqnn
M(t)=\int_0^t \phi(X(s))\cdot dW_s\,,
\eeqnn
$\langle M\rangle$ is its bracket
\beqnn
\langle M\rangle(t)=\int_0^t \vert \phi(X(s))\vert^2\, ds\,,
\eeqnn
and $\phi(x)=\sigma^{-1}(x)(c(x)-b(x))$.

In the next discussion we use the expression \lq\lq Brownian motion with covariance $\sigma^2$ and drift $c$\rq\rq\ to denote
any process whose law is the same as the law of $(\sigma W_t+ct\,;\,t\geq 0)$.

If we choose $\sigma$ and $c$ constant and $b=0$ above, then $Y$ is a
Brownian motion with covariance $\sigma^2$
and drift $c$ and $(X,M)$ is a centered Brownian motion (in dimension $d+1$)
whose covariance satisfies the following:
the covariance of $X$ is $\sigma^2$; $E[M(t) X(t)]=ct$.
Thus the Girsanov formula then has the following corollary:  let $Y$ be a Brownian motion
with covariance matrix $\sigma^2$ and drift $c$, then
\beqn \label{eq:malliavin}
E[F(Y([0,t]))]=E[F(X([0,t]))\,e^{M(t)-\frac{1}{2}\langle M\rangle (t)}]\,,
\eeqn
for any time $t$,
for any bounded continuous function $F$ on $C([0,t],\R^d)$ and
for any random process $M$ such that $(X,M)$ is a centered Brownian motion, $X$ has covariance $\sigma^2$
and the covariance of $X$ and $M$ is $E[M(t) X(t)]=ct$.

Applying the Girsanov formula to the processes $X^\o$ and
$X^{\l,\,\o}$, we get that, for any $\o$, \beqn\label{eq:gir}
E[F(X_0^{\l,\,\o}([0,t]))]=E[F(X_0^{\o}([0,t]))\, e^{\l\bB (t)
-\frac {\l^2}2 \langle\bB\rangle(t)}]\,, \eeqn where $\bB$ is the
martingale
$$ \bB (t)=\int_0^t \sigma^\o(X_0^\o(s))\, e_1 \cdot dW_s$$
and $\langle\bB\rangle$ is its bracket
$$\langle\bB\rangle(t)=\int_0^t \vert \sigma^\o(X_0^\o(s))\, e_1 \vert^2 \, ds \,.$$

In particular, in the range  $\l^2t=\alpha$, we have
\beqnn
E[F(X_0^{\l,\,\o}([0,t]))]=E[F(X_0^{\o}([0,t]))\,
e^{\l\bB (\frac \alpha{\l^2}) -\frac {\l^2}2 \langle\bB\rangle(\frac \alpha{\l^2})}]\,.
\eeqnn

We shall need the following easy statement:

\begin{lemma} \label{lm:ap1}
For all $\alpha\geq 1$ and $\mu>1$, we have
\beqn\label{eq:apriori1}
\limsup_{\l\rightarrow 0\,;\, t\rightarrow +\infty \,;\, \l^2t=\alpha}
E\left[\, e^{\mu\l\bB (\frac \alpha{\l^2}) -\frac {\mu\l^2}2 \langle\bB\rangle(\frac \alpha{\l^2})}\right]
<\infty\,,
\eeqn
uniformly in $\o$.
\end{lemma}

\medskip\noindent
{\it Proof of (\ref{eq:apriori1})}

Assumption 3 on the ellipticity of $a^\o$ implies
that $\langle\bB\rangle(t)\leq\kappa^{-1}t$. Therefore
\beqnn
&&E\left[\, e^{\mu\l\bB (\frac \alpha{\l^2}) -\frac {\mu\l^2}2 \langle\bB\rangle(\frac \alpha{\l^2})}\right]\\
&\leq& E\left[\, e^{\mu\l\bB (\frac \alpha{\l^2}) -\frac {\mu^2\l^2}2 \langle\bB\rangle(\frac \alpha{\l^2})}\right]
e^{(\mu-1)\frac {\mu\l^2}2 \frac \alpha{\kappa\l^2}} \\
&=&e^{(\mu-1)\frac {\mu}2 \frac \alpha{\kappa}}\,,
\eeqnn
and (\ref{eq:apriori1}) is thus proved. \qed

Next we apply a (joint) invariance principle for the process
$(X_0^{\o,\,\eps}, \bB^\eps)$
where $X_0^{\o,\,\eps}(t)=\eps X_0^\o(t/\eps^2)$ and $\bB^\eps(t)=\eps\bB(t/\eps^2)$.

Let us recall some of the ideas of \cite{kn:KV} and \cite{kn:DFGW}.

The process of the environment seen from the particle: $(\o(t)=X_0^\o(t).\o)\,;\,t\geq 0)$ is
a Markov process under the annealed law, with values in $\O$.
It is not difficult to check that the measure $\Q$
is invariant, ergodic and reversible for this process (Recall that $V=0$ for now!).

Given the state of the environment at times $0$ and $t$, say $\o(0)$ and $\o(t)$, one retrieves the position of
the particle itself by solving the equation $z.\o(0)=\o(t)$. Note that Assumptions 1 (either $V$ or $a$ is not constant)
and 4 (independence property) imply that there cannot be more than one solution.
It also follows from the equality $(X_0^\o(t)-X_0^\o(s)).\o(s)=\o(t)$ that $X_0^\o(t)$
is an antisymmetric additive functional of the process $\o(\cdot)$.
(Antisymmetric means that reversing time amounts to
changing the sign of $X_0^\o$.)

The process $\bB$ is also an additive functional of $\o(\cdot)$ since it can be written
as the difference
\beqnn \bB(t)-\bB(s)&=&\tr{e_1}\, (X_0^\o(t)-X_0^\o(s))-\int_s^t \tr{e_1}\, b^\o(X_0^\o(u))\,du\\
&=& \tr{e_1}\, (X_0^\o(t)-X_0^\o(s))-\int_s^t \tr{e_1}\, b(\o(u))\,du\,,
\eeqnn
see (\ref{eq:sde1}).

Sufficient conditions for invariance principles for additive functionals of reversible Markov
processes that can be applied to $X_0^\o$ or $\bB$ are given in \cite{kn:KV} and \cite{kn:DFGW}. They
yield a joint invariance principle for $(X_0^{\o,\,\eps}, \bB^\eps)$.

In order to compute the covariance matrix note that
$X_0^\o(t)$ is antisymmetric whereas
$\int_0^t \tr{e_1}\, b(\o(s))\,ds$ is a symmetric functional of the environment. Thus they are orthogonal
under the annealed measure. Therefore
$$\int E[\bB(t)\, X_0^\o(t)]\, d\Q=
\int E[( \tr{e_1}\, X_0^\o(t))\, X_0^\o(t)]\, d\Q\,, $$
and  thus the asymptotic covariance of $X_0^{\o,\,\eps}(1)$ and $\bB^\eps(1)$
coincides with the asymptotic covariance of $X_0^{\o,\,\eps}(1)$
and $\tr{e_1}\, X_0^{\o,\,\eps}(1)$
and equals $\Sigma e_1$.

Applying first the Girsanov formula and then the invariance principle, we get that,
for any $\alpha>0$, as $\l\to 0$ and $t\to \infty$ with $\l^2 t = \alpha$, we have
\begin{eqnarray}\label{eq:girs1}
&&\int E\left[F\left(\frac \l \alpha X_0^{\l,\,\o}(\frac \alpha{\l^2}s)\,;\, 0\leq s\leq 1 \right)
\, \right]d\Q\cr
&=& \int E\left[F\left(\frac \l \alpha X_0^\o(\frac \alpha{\l^2}s)\,;\, 0\leq s\leq 1 \right)
e^{\l\bB (\frac \alpha{\l^2}) -\frac {\l^2}2 \langle\bB\rangle(\frac \alpha{\l^2})}\right] \, d\Q\cr
&&\rightarrow E\left[F\left(\frac 1 {\sqrt \alpha} N(s)\,;\, 0\leq s\leq 1 \right)
e^{\sqrt{\alpha}Z(1)-\frac \alpha2 E(Z(1)^2)}\right]\,,
\end{eqnarray}
where $F$ is a bounded continuous functional on $C([0,1],\R^d)$ and
$(N,Z)$ is a centered Brownian motion of dimension $d+1$
with $N$ having $\Sigma$ as covariance matrix, and $E[Z(t)N(t)]=\Sigma e_1 t$.
We refer to Lemma \ref{lm:ap1} for the full justification of the passing to the
limit in (\ref{eq:girs1}).
Using now formula (\ref{eq:malliavin}), we have:
\beqnn
E\left[F\left(\frac 1 {\sqrt \alpha} N(s)\,;\, 0\leq s\leq 1 \right)
e^{\sqrt{\alpha}Z(1)-\frac \alpha2 E(Z(1)^2)}\right]
= E\left[F\left(\frac 1 {\sqrt \alpha}N(s)+\Sigma\, e_1 s\,;\, 0\leq s\leq 1 \right)\right]\,.
\eeqnn
We thus conclude that
\beqnn &&\E_0^\l\left[F\left(\frac \l \alpha X(\frac \alpha{\l^2}s)\,;\, 0\leq s\leq 1  \right)\right]
=\int E\left[F\left(\frac \l \alpha X_0^{\l,\,\o}(\frac \alpha{\l^2}s)\,;\, 0\leq s\leq 1 \right)
\, \right]d\Q\cr
&&\rightarrow E\left[F\left(\frac 1 {\sqrt \alpha}N(s)+\Sigma\, e_1 s\,;\, 0\leq s\leq 1 \right)\right]\,, \eeqnn
i.e., when  $\l^2t=\alpha$, the law of
$((\l t)^{-1} X(ts)\,;\, 0\leq s\leq 1 )$ under $\P_0^\l$ converges to the law of
$(\frac 1 {\sqrt \alpha}N(s)+\Sigma\, e_1 s\,;\, 0\leq s\leq 1 )$.

To finish the proof of the Proposition, we need a priori bounds on the moments of
$\vert X(t)\vert$ under $\P^\l_0$.  We shall prove in Lemma \ref{lm:ext} that
\beqnn
\limsup_{\l\rightarrow 0\,;\, t\rightarrow +\infty \,;\, \l^2t=\alpha}
E_0^{\l,\,\o}\left[\,\frac {\max\limits_{s \le t}\vert X(s)\vert^p}{\l^p t^p}\right]
<\infty\,,
\eeqnn
uniformly in $\o$ and for all $p\geq 1$ and all $\alpha\geq 1$.
Therefore
\beqn\label{eq:apriori}
\limsup_{\l\rightarrow 0\,;\, t\rightarrow +\infty \,;\, \l^2t=\alpha}
\E_0^\l\left[\,\frac {\max\limits_{s \le t}\vert X(s)\vert^p}{\l^p t^p}\right]
<\infty\,,
\eeqn
for all $p\geq 1$ and all $\alpha\geq 1$
and we observe that (\ref{eq:apriori}) together with the convergence of the law of
$((\l t)^{-1} X(ts)\,;\, 0\leq s\leq 1 )$ under $\P_0^\l$ to the law of
$(\frac 1 {\sqrt \alpha}N(s)+\Sigma\, e_1 s\,;\, 0\leq s\leq 1 )$
implies (\ref{eq:***}) and (\ref{eq:max****}).
Indeed we have
$$\frac {\max\limits_{s \le t}\vert X(s)\vert^2}{\l^2 t^2}
=\max\limits_{s \le 1}\vert \frac \l\a X(\frac \alpha{\l^2}s)\vert^2$$
and therefore, with the notation above,
$$\E_0^\l\left[\,\frac {\max\limits_{s \le t}\vert X(s)\vert^2}{\l^2 t^2}\right]
\rightarrow
E\left[\,\max\limits_{s \le 1}\vert \frac 1 {\sqrt \alpha}N(s)+\Sigma\, e_1 s\vert^2\right]\,,$$
where we  used the function $F(w(s)\,;\, 0\leq s\leq 1)=\max\limits_{s \le 1}\vert w(s)\vert^2$ and
inequality (\ref{eq:apriori}) to justify the passing to the limit.
Finally it is easy to check that
$$\sup_{\alpha\geq 1}
E\left[\,\max\limits_{s \le 1}\vert \frac 1 {\sqrt \alpha}N(s)+\Sigma\, e_1 s\vert^2\right]<\infty\,.$$
This last line ends the justification of (\ref{eq:max****}).
Equation (\ref{eq:***}) is proved the same way using now the function $F(w(s)\,;\, 0\leq s\leq 1)=w(1)$.

The proof of Proposition \ref{lm:LR} is now complete in the case $V=0$. \qed

\medskip\noindent
{\it Proof of Proposition \ref{lm:LR} in the case $V\not=0$}

We do not assume anymore that $V=0$.

Define $Y^\o$ and $ Y^{\l,\,\o}$ to be the solutions of the stochastic differential equations

\beqn \label{eq:sde11}
dY^\o(t)=e^{-2V^\o(Y^\o(t))} b^\o(Y^\o(t))\,dt+e^{-V^\o(Y^\o(t))} \sigma^\o(Y^\o(t))\,dW_t \,;\, Y^\o(0)=0\,,
\eeqn
and
\begin{eqnarray}  \label{eq:sde21}
dY^{\l,\,\o}(t)=&&e^{-2V^\o(Y^{\l,\,\o}(t))} b^\o(Y^{\l,\,\o}(t))\,dt+
\,e^{-2V^\o(Y^{\l,\,\o}(t))} a^\o(Y^{\l,\,\o}(t))\hl\,dt\cr &&
+e^{-V^\o(Y^{\l,\,\o}(t))} \sigma^\o(Y^{\l,\,\o}(t))\,dW_t \,;\, Y^{\l,\,\o}(0)=0\,,
\end{eqnarray}

so that the generators of $Y^\o$ and $ Y^{\l,\,\o}$ are the operators
\beqnn \MM^\o f(x)=\frac 1 2 \, \mathrm{div}(e^{-2V^\o}\,a^\o\,\nabla f)(x)\,,
\eeqnn
and
\beqnn
 \MM^{\l,\,\o} f(x)
=
\frac 1 2 \, e^{\hl\cdot x}\mathrm{div}(e^{-2V^{\l,\,\o}}\,a^\o\,\nabla f)(x)\,,
\eeqnn
where $V^{\l,\,\o}(x)=V^\o(x)-\tr{\hl}\, x$.

Note that these operators are of the same form as $\LL^\o$ and $\LL^{\l,\,\o}$ with
$V^\o$ being replaced by $0$ and $a^\o$ being replaced by $\exp(-2V^\o)a^\o$.
Thus we may apply the results obtained in the special case of a vanishing potential to the processes
$Y^\o$ and $Y^{\l,\,\o}$, in particular $Y^\o$ satisfies the invariance principle with some asymptotic
diffusivity $\Sigma^Y$ and $Y^{\l,\,\o}$ satisfies:

\beqn \label{eq:***'}
\lim_{\l\rightarrow 0\,;\, t\rightarrow +\infty \,;\, \l^2t=\alpha}
\int_\Omega E\left[\,\frac {Y^{\l,\,\o}(t)}{\l t}\right]d\Q
= \Sigma^Y\, e_1\,,
\eeqn
and
\beqn
\label{eq:max****'}
\sup_{\a\geq 1}\limsup_{\l\rightarrow 0\,;\, t\rightarrow +\infty \,;\, \l^2t=\alpha}
\int_\Omega
E\bigg[\,\frac {\max\limits_{s \le t}\vert Y^{\l,\,\o}(s)\vert^2}{\l^2 t^2}\bigg]d\Q
<\infty\,.
\eeqn

Fix $\o$ and
set $$A^\o(t)=\int_0^t e^{-2V^\o(Y^\o(s))}\, ds
=\int_0^t e^{-2V(Y^\o(s).\o)}\, ds\,.$$
Then $X_0^\o$ has the same law as the time changed process $Y^\o((A^\o)^{-1})$.
Similarly, if we let $$ A^{\l,\,\o}(t)=\int_0^t e^{-2V^\o(Y^{\l,\,\o}(s))}\, ds\,,$$
then $X_0^{\l,\,\o}$ has the same law as the time changed process $Y^{\l,\,\o}((A^{\l,\,\o})^{-1})$.

>From Assumption 3, we know that $V$ is bounded and therefore \beqn\label{eq:A} A^{\l,\,\o}(t)\geq \cinf t\,,\eeqn
for some constant $\cinf\leq 1$.
Observe that  (\ref{eq:max****}) immediately follows from (\ref{eq:A}) and
 (\ref{eq:max****'}).

\medskip\noindent
{\it Proof of (\ref{eq:***})}

Let
$$\g=\int d\Q(\o)\, e^{-2V(\o)}\,.$$
The ergodic theorem for the process $Y^\o.\o$ implies that $A^\o(t)/t$ almost surely converges to $\g$. We need
a similar  statement for $A^{\l,\,\o}$:

\begin{lemma} \label{lm:Alo}
For fixed $\alpha>0$ and any positive $\eta$, we have $\Q$--a.s.
\beqn \label{eq:Alo}
P\left[\Big\vert \frac{A^{\l,\,\o}(t)}t-\g\Big\vert\geq\eta\right]\rightarrow 0\,,
\eeqn
for $\l\rightarrow 0$, $t\rightarrow +\infty$, $\l^2t=\alpha$.
\end{lemma}

\noindent
{\it Proof of Lemma \ref{lm:Alo}}\\
Use the Girsanov formula (\ref{eq:gir}) to see that
\beqnn
P\left[\Big\vert \frac{A^{\l,\,\o}(t)}t-\g\Big\vert\geq\eta\right]
=E\left[\1(\Big\vert \frac{A^\o(t)}t-\g\Big\vert\geq\eta) e^{\l \bB (t) -\frac{\l^2}2 \langle\bB\rangle(t)}\right]\,.
\eeqnn
The convergence in (\ref{eq:Alo}) then follows from H\"older's inequality, the bound in Lemma \ref{lm:ap1} and
the fact that $P[\vert \frac{A^\o(t)}t-\g\vert\geq\eta]\rightarrow 0$. \qed

Back to the proof of (\ref{eq:***}), we start with the equality
\beqnn
\E_0^\l\left[\,\frac {X(t)}{\l t}\right]
=
\int_\Omega E\left[\,\frac {Y^{\l,\,\o}((A^{\l,\,\o})^{-1}(t))}{\l t}\right]\, d\Q\,.
\eeqnn
We have
\begin{eqnarray*}
&&\Big\vert \int_\Omega E\left[\,\frac { Y^{\l,\,\o}((A^{\l,\,\o})^{-1}(t)) } {\l t}\right]d\Q
- \int_\Omega E\left[\,\frac {Y^{\l,\,\o}(t/\g)} {\l t}\right]d\Q \Big\vert\\[4mm]
&\leq&
\int_\Omega E\Bigg[ \frac{ \max\limits_{\vert s -t/\g\vert\leq \eta t}\vert Y^{\l,\,\o}(s)-Y^{\l,\,\o}(t/\g)\vert }{\l t}\Bigg]d\Q\\[3mm]
&+&\int_\Omega  E\Bigg[ \frac{ \max\limits_{s\leq  t/\cinf}\vert Y^{\l,\,\o}(s)\vert} {\l t}
\1(\vert (A^{\l,\,\o})^{-1}(t)-t/\g\vert\geq\eta t)\Bigg]d\Q \\[2mm]
&=& I+II\,.
\end{eqnarray*}

($\cinf$ is the same constant as in (\ref{eq:A}). Note that $\cinf\leq \g$.)

By the Markov property, we have
\beqnn
I\leq
2\int_\Omega
E\Bigg[ \frac{ \max\limits_{s \leq 2\eta t}\vert Y^{\l,\,\o}(s)\vert}{\l t}\Bigg]d\Q\,.
\eeqnn
Now (\ref{eq:max****'}) implies that, if we let $\l$ tend to $0$
and then  $\eta$ tend to $0$, then
the contribution of I vanishes.
H\"older's inequality, (\ref{eq:max****'}) again and Lemma \ref{lm:Alo} imply that, for any $\eta>0$,
then II also converges to $0$ as $t$ tends to $+\infty$.
We conclude that
\beqnn
\lim_{\l\rightarrow 0\,;\, t\rightarrow +\infty \,;\, \l^2t=\alpha}
\int_\Omega \left(E\left[\,\frac {Y^{\l,\,\o}((A^{\l,\,\o})^{-1}(t))}{\l t}\right]
- E\left[\,\frac {Y^{\l,\,\o}(t/\g)}{\l t}\right] \right)d\Q=0\,,
\eeqnn
and, using (\ref{eq:***'}),
\beqnn
\lim_{\l\rightarrow 0\,;\, t\rightarrow +\infty \,;\, \l^2t=\alpha} \E_0^\l\left[\,\frac {X(t)}{\l t}\right]
=\frac 1\g\, \Sigma^Y\, e_1\,.\eeqnn

The last piece of information missing is the equality $\Sigma^Y=\g\Sigma$. It comes as follows:
since $\frac 1{\sqrt{t}} Y^\o(t)$ converges in law to a Gaussian random variable with covariance
$\Sigma^Y$ and since $\frac 1 t A^\o(t)$ almost surely converges to $\g$ then
$\frac 1{\sqrt{t}} X_0^\o(t)=\frac 1{\sqrt{t}} Y^\o((A^\o)^{-1}(t))$ converges in law to a Gaussian random variable with covariance
$\Sigma=\frac 1\g \Sigma^Y$.\qed


\vskip 1cm
\section{A priori estimates}
\setcounter{equation}{0}
\label{sec:IV}

In this section, we prove some a priori estimates on exit times
that quantify the fact that the process $X^{\l,\,\o}$
is transient in the direction $e_1$.

For a given realization of the environment, the local drift of the process $X^\o$
equals $b^\o(X^\o(t))$. Its mean under the annealed law vanishes. The drift of
$X^{\l,\,\o}$ has an extra $a^\o(X^{\l,\,\o}(t))\hl$ term. Since, by Assumption
3 (uniform ellipticity), we have $\tr{e_1}\, a^\o(x) e_1\geq \kappa \vert e_1\vert^2\not=0$
for any $\o$ and $x$, one would expect $X^{\l,\,\o}$ to be transient in the direction $e_1$,
and this turns out to be the case, but we also need more quantitative statements on the
tendency of the diffusion to go in the direction $e_1$.

Roughly speaking, we may think
of $\tr{e_1}\, X^{\l,\,\o}(t)$ as the sum of a centered term of order $\sqrt{t}$ and a
drift term of order $\l t$. Thus the shortest scale on which we may hope the drift term to
dominate is $\l^2 t\geq 1$ or, in terms of space scale, $\l L \geq 1$.
Up to the value of the constants $\kappa_1$, $\kappa_2$, $c$ and $C$, our estimates
in Lemma \ref{lm:leftexit} and Lemma \ref{lm:rightexit} are therefore optimal.

In the following Lemmata, inequalities (\ref{eq:leftexit}) and (\ref{eq:rightexit}) hold true for any environment $\omega$
satisfying Assumptions 2 and 3. Assumptions 1 and 4 are not relevant in this section.

We use $T_L=\inf\{t\,: \, \tr{e_1}\,X(t)=L\}$ to denote the hitting time of the hyperplane $\{x: \tr{e_1}x =L\}$, $L \in \R$.

\begin{lemma} \label{lm:leftexit}
There exists constants $c>0$ and $\kappa_1>0$
that depend on the dimension, the ellipticity constant $\kappa$
and the $L^\infty$ bound on $V$ such that for all $L$, $\l\leq 1$ and for any environment $\o$,
\beqn\label{eq:leftexit}
P_0^{\l,\,\o}[\,T_{-L}<\infty]\leq c e^{-\kappa_1 \l L}\,.
\eeqn
\end{lemma}

\begin{lemma} \label{lm:rightexit}
There exists constants $C$ and $\kappa_2>0$ that depend on the dimension, the ellipticity constant $\kappa$
and the $L^\infty$ bound on $V$ such that for all $L$, $\l\leq 1$, $t\geq 0$ and
for any environment $\o$,
\beqn\label{eq:rightexit}
P_0^{\l,\,\o}[\,T_L\geq t]\leq Ce^{-\kappa_2\l^2 t +\l L}\,.
\eeqn
\end{lemma}

\medskip\noindent
{\it Proof of Lemma \ref{lm:leftexit}}.\\
Let $u(x)= u_{\l, \o, L}(x): = P_x^{\l,\,\o}[\,T_L< T_{-L}]$. Then, $u$ solves the elliptic boundary value problem
\begin{equation}\label{aux_pde}
 \LL^{\l,\,\o} u(x) =0, \quad -L \leq \tr{e_1}\,x \leq L\, ,
\end{equation}
where the generator $\LL^{\l,\,\o}$ was defined in (\ref{generator}), with boundary values
\begin{equation}\label{aux_bc}
u(x) = 0 \, \hbox{ if } \tr{e_1}\,x = -L \hbox{ and } u(x) = 1 \, \hbox{ if }\tr{e_1}\,x = L\, .
\end{equation}
We first need the following Lemma, whose proof is deferred.
\begin{lemma} \label{lm:withoutl}
Assume $\l =1$ and $u$ is the solution of the boundary value problem above. Then, under Assumption 3, there exists $L_0 > 0$ depending only on $\kappa$ and on the dimension, such that for all $L\geq L_0$ and all $\o$,
\begin{equation}\label{woutl}
u(x) \geq \frac{2}{3} \quad \hbox{ for all } x \hbox{ with } \tr{e_1}\,x = 0\, .
\end{equation}
\end{lemma}
Let $\widetilde{X}(t) = \l X\left( \frac{t}{\l^2}\right), t \geq 0$. Then, $(\widetilde{X}(t))_{t \geq 0}$ is a Markov process with generator
$$
\widetilde{\LL}^{\l,\,\o} f(x)=
\frac 1 2 e^{2V^\o(x/\l)}\, \mathrm{div}(e^{-2V^\o(x/\l)}\,a^\o(x/\l)\,\nabla f)(x)+\,a^\o(x/\l)e_1\,\cdot \nabla f(x)
$$
(The advantage of $(\widetilde{X}(t))_{t \geq 0}$ is that we scaled away the drift and we will be able to apply Lemma \ref{lm:withoutl}). Let us consider exit times for $(\widetilde{X}(t))_{t \geq 0}$.
$$
\widetilde{T}_L=\inf\{t\,: \, \tr{e_1}\,\widetilde{X}(t)=L\}, \quad
\widetilde{T}_{-L}=\inf\{t\,: \, \tr{e_1}\,\widetilde{X}(t)=-L\},
$$
and
$$
\widetilde{T}_{\pm L}=\inf\{t\,: \, \vert \tr{e_1}\,\widetilde{X}(t)\vert =L\}\, .
$$
Then,
$$
P_0^{\l,\,\o}[\,T_{-L}<\infty] = P_0^{\l,\,\o}[\,\widetilde{T}_{-\l L}<\infty]
$$

Hence, it suffices to show that for all $\l, L, \o$, we have $P_0^{\l,\,\o}[\,\widetilde{T}_{-\l L}<\infty] \leq c e^{-\kappa_1 \l L}$ which in turn is equivalent to
\begin{equation}\label{boundwoutl}
P_0^{\l,\,\o}[\,\widetilde{T}_{-L}<\infty] \leq c e^{-\kappa_1 L} \, ,
\end{equation}
for all $\l, L, \o$.
Let $L_0$ be as in Lemma \ref{lm:withoutl}. We consider the embedded random walk defined as follows.
Let $t_1 : = \widetilde{T}_{\pm L_0}$, $t_{i+1}= \inf\{t \geq t_i: \vert \tr{e_1}\,(\widetilde{X}(t)- X(t_i))\vert =L_0\}$ and $S_i = \widetilde{X}_{t_i}, i=1,2, \ldots $, (and $S_0= 0$). Due to Lemma \ref{lm:withoutl} and the strong Markov property of $(\widetilde{X}(t))_{t \geq 0}$, we have
$$
P_0^{\l,\,\o}[S_{i+1} = S_i +L_0] \geq \frac{2}{3}\, .
$$
Hence we can couple $(S_i)_{i=0, 1,2, \ldots}$ with a standard random walk with drift
$(\bar{S}_i)_{i=0, 1,2, \ldots}$ on $L_0 \Z $ with iid increments, satisfying
$P[\bar{S}_{i+1} = \bar{S}_i +L_0] = \frac{2}{3} = 1-P[\bar{S}_{i+1} = \bar{S}_i -L_0]$
and the coupling is such that $S_i\geq \bar{S}_i$ for all $i$.
Explicit calculation yields
$$
P[\bar{S}_i > -m L_0, \, \forall i] = 1 - \frac{1}{2^m}\, ,
$$
and we obtain
$$
P_0^{\l,\,\o}[S_i > -m L_0, \, \forall i] \geq 1 - \frac{1}{2^m}\, ,
$$
which implies (\ref{boundwoutl}).\qed

\medskip\noindent
{\it Proof of Lemma \ref{lm:rightexit}}.\\
Inequality (\ref{eq:rightexit}) is equivalent to the following statement in terms of $(\widetilde{X}(t))_{t \geq 0}$:
There exist constants $C$ and $\kappa_2>0$ that depend on the dimension, the ellipticity constant $\kappa$
and the $L^\infty$ bound on $V$ such that for all $L$ and $t$
and for any environment $\o$,
\begin{equation}\label{rightexittild}
P_0^{\l,\,\o}[\,\widetilde{T}_L\geq t]\leq Ce^{-\kappa_2 t + L}\,.
\end{equation}

Indeed, note that $\widetilde{T}_{\l L} = \l^2 T_L$.

Take $L_0$ as in Lemma \ref{lm:withoutl}. Then Aronson's estimate (see \cite{kn:aronson}) yields
 that
for all $x \in \Pi_{L_0}: = \{x: -L_0 \leq \tr{e_1}\,x \leq L_0\}$,
\begin{equation}\label{aroncons}
P_x^{\l,\,\o}[\vert e_1\cdot\widetilde{X}(1)\vert > L_0] \geq \gamma > 0\, ,
\end{equation}
where $\gamma$ only depends on $L_0$, $\kappa$ and the dimension.
Indeed, according to \cite{kn:aronson}, the fundamental solution ${\cal P}(x,y,t)$ of the parabolic operator
$$
\partial_t-\mathrm{div}(a(x)\nabla)-a(x)e_1\cdot \nabla
$$
(or, equivalently, transition probability density of the corresponding diffusion) satisfies the lower bound
$$
{\cal P}(x,y,1)\ge k\exp(-K|x-y|^2)
$$
with constants $k>0$ and $K>0$ which only depend on $\kappa$ and $d$.
The estimate (\ref{aroncons}) is an immediate consequence of this lower bound.
Due to the Markov property of $(\widetilde{X}(t))_{t \geq 0}$,  estimate (\ref{aroncons}) implies that for some constant $g_0$ (which depends only on $L_0$, $\kappa$ and the dimension),
\begin{equation}\label{fastexit}
P_x^{\l,\,\o}[\widetilde{T}_{{\pm L}_0} \geq g_0] \leq \frac{1}{15}\, ,
\end{equation}
for all $x \in \Pi_{L_0}$. Define the stopping time $\widetilde{t}$ as follows:
$$
\widetilde{t}: = \inf\{t \geq 0: \vert \tr{e_1}\,(\widetilde{X}(t) - \widetilde{X}(0))\vert = L_0\}
\wedge g_0\, .
$$
Combining (\ref{fastexit}) with Lemma \ref{lm:withoutl}, we obtain
$$
P_x^{\l,\,\o}\left[\tr{e_1}\,\widetilde{X}(\widetilde{t}) = \tr{e_1}\,x + L_0\right]\geq \frac{3}{5}\,
$$
(with probability at least $\frac{14}{15}$, the layer $\{y: \vert \tr{e_1}\,y - \tr{e_1}\,x \vert < L_0\}$ has been left by time $g_0$, and with probability at least $\frac{2}{3}$, the exit happens at $\{y: \tr{e_1}\,y = \tr{e_1}\,x + L_0\}$).
We consider the embedded random walk defined as follows.
Let $t_1 : = \widetilde{t}$, $t_{i+1}= \widetilde{t}\circ \theta_{t_i}+ t_i$, $S_i = \tr{e_1}\widetilde{X}(t_i), i=1,2, \ldots $ and $S_0 =0$.
(Here $\theta$ denotes the shift operator on path space). Hence we can couple $(S_i)_{i=0, 1,2, \ldots}$ with a standard random walk with drift
$(\bar{S}_i)_{i=0, 1,2, \ldots}$ on $L_0 \Z $ with i.i.d. increments, satisfying
$P[\bar{S}_{i+1} = \bar{S}_i +L_0] = \frac{3}{5} = 1-P[\bar{S}_{i+1} = \bar{S}_i -L_0]$
and the coupling is such that $S_i\geq \bar{S}_i$ for all $i$. It is straightforward to check that there are constants $\kappa_3 > 0$ and $c> 0$ such that for all $n \in \N$ and $L \in \R$,
$$
P[\bar{S}_n \leq L L_0]\leq ce^{- \kappa_3 n + L}\, .
$$
We conclude, by comparison, that we have for all $L \in \R$ and $n \in \N$
$$
P_0^{\l,\,\o}[\widetilde{T}_{L L_0}\geq n g_0] \leq P_0^{\l,\,\o}[S_n \leq L L_0] \leq P[\bar{S}_n \leq L L_0]\leq ce^{- \kappa_3 n+L}\, ,
$$
and this implies (\ref{rightexittild}).\qed

\medskip\noindent
{\it  Proof of Lemma \ref{lm:withoutl}.}\\
Without loss of generality we can assume that $V=0$.
Indeed, multiplying (\ref{aux_pde}) by $\exp(-2V^\omega(x))$ and
denoting $a(x)= \exp(-2V^\omega(x))a^\omega(x)$ one can rewrite  equations
(\ref{aux_pde})--(\ref{aux_bc}) in the form
\begin{equation}\label{eq0}
\mathrm{div}(a(x)\nabla u)+a(x)e_1\cdot \nabla u=0\qquad\hbox{in }\Pi,
\end{equation}
$$
u(-L,z)=0,\qquad u(L,z)=1.
$$
with $x =(x_1, z)$ and  $\Pi=\Pi_L=\{x\in \mathbb R^d\,:\,-L\le x_1\le L\}$. Under  Assumption 3 the matrix $a(x)$ is symmetric and
satisfies the
following elliptic estimates
\begin{equation}\label{uniell}
\kappa |y|^2 \le a(x)y\cdot y\le \kappa^{-1} |y|^2,\qquad \forall x,\,y \in\mathbb R^d.
\end{equation}
Therefore, it suffices to prove the following statement.
\begin{lemma}\label{l_main}
There is $L_0=L_0(\kappa,d)$ such that for all $L\ge L_0$, we have
$$
u(0,z)\ge 2/3,\quad z \in\mathbb R^{d-1}.
$$
\end{lemma}
{\it Proof of Lemma \ref{l_main}}\\
Suppose that $u(0,z)<2/3$ for some $z\in\mathbb R^{d-1}$ and some $L$.
Without loss of generality we assume that $z=0$. The function
$1-u(x)$ is a non-negative solution of the equation
$$
\mathrm{div}(a(x)\nabla (1-u))+a(x)e_1\cdot \nabla(1-u)=0,
$$
therefore, by the Harnack inequality (see, for instance, \cite{kn:GT}) there is a constant $C=C(\kappa,d)$ such that
$$
1-u(0,z)\ge C(\kappa,d)(1-u(0))\ge \frac{1}{3}C(\kappa,d),
\qquad z \in [-1,1]^{d-1}.
$$
For all $z\in [-1,1]^{d-1}$ this implies the estimate
\begin{eqnarray}
\frac{1}{9}C^2(\kappa,d) &\le& \bigg(\int\limits_0^L \frac{\partial(1-u(x_1,z))}{\partial x_1}\,dx_1\bigg)^2\le
L\int\limits_0^L \Big(\frac{\partial u(x_1,z)}{\partial x_1}\Big)^2 \,dx_1 \cr
&\le&
L\int\limits_0^N |\nabla u(x_1,z)|^2 \,dx_1\le L\int\limits_{-L}^L e^{x_1}|\nabla u(x_1,z)|^2 \,dx_1,
\end{eqnarray}
where we used Jensen's inequality for the second inequality.
Integrating over $[-1,1]^{d-1}$, we obtain
\begin{equation}\label{lower}
\int\limits_{-L}^L dx_1\int\limits_{[-1,1]^{d-1}}e^{x_1} |\nabla u(x_1,z)|^2 \,dz\ge \frac{1}{L}C_1(\kappa,d),
\end{equation}
with a constant $C_1(\kappa,d)>0$ which only depends on $\kappa$ and $d$.
For arbitrary sets $G\subseteq \Pi$ and functions $v(x)$, we write
$\EE(v,G)$ for the energy of $v$ on $G$:
$$
\EE(v,G)=\int_G e^{x_1}a(x)\nabla v(x)\cdot\nabla v(x)\,dx.
$$
Hence, (\ref{lower}) can be rewritten, setting $G_0 = [-L, L]\times [-1,1]^{d-1}$, as follows:
\begin{equation}\label{lowereng}
\EE(u, G_0) \ge \frac{1}{L}C_1(\kappa,d)\, .
\end{equation}
(The value of $C_1(\kappa,d)$ changes from equation (\ref{lower}) to equation (\ref{lowereng}).)
We will give an upper bound for $\EE(u, G_0)$ which contradicts (\ref{lowereng}) when $L$ is too large.
We introduce the following subsets of $\Pi$.
$$
B_0=\big[-e^{L/d},e^{L/d}\big]^{d-1},\qquad
\Pi_0=[-L,L]\times B_0,
$$
$$
B_j=2^{-j}B_0=2^{-j}\big[-e^{L/d}, e^{L/d}\big]^{d-1},\qquad \Pi_j=[-L,L]\times B_j,
$$
Also, denote
$$
\bar u(x)=\bar u(x_1)=\left\{\begin{array}{ll}
1,\qquad &\hbox{if } x_1\ge -L+1,\\[2mm]
x+L,\quad&\hbox{if } -L\le x_1\le-L+1\end{array}\right.
$$
and, writing again $x = (x_1, z)$,
$$
v_j(x)=\left\{\begin{array}{ll}
(1-d_j(z))u(x)+d_j(z)\bar u(x),\qquad &\hbox{if } d_j(z)\le 1,\\[2mm]
\bar u(x),\quad&\hbox{otherwise},\end{array}\right.
$$
where
$$
d_j(z)=2^{j+1}e^{-L/d}\,\mathrm{dist}_{\mathbb R^{d-1}}(z,\partial B_j),\qquad j=0,1,\dots,2d.
$$
Note that $v_j=\bar u$ in the domain $\Pi_{j+1}$.
Since $v_0=u$ on $\partial \Pi_0$, we know that
$$
\EE(u,\Pi_0)\le \EE(v_0,\Pi_0)\, .
$$
Indeed, the functional $\{\EE(v,\Pi_0):v\in H^1(\Pi_0),\,v\big|_{\partial\Pi_0}=u\big|_{\partial\Pi_0}\}$ attains its minimum at the unique solution of the equation
$$
\mathrm{div}(a(x)\nabla v)+a(x)e_1\nabla v=0\quad\hbox{in }\Pi_0,\qquad
v\big|_{\partial\Pi_0}=u\big|_{\partial\Pi_0}.
$$
Clearly, the function $u$ solves this equation, and the required inequality follows.

For $x\in\Pi_0\setminus\Pi_1$ we have
$$
e^{x_1}a(x)\nabla v_0(x)\cdot\nabla v_0(x)
$$
$$
=
e^{x_1}a(x) \{(1-d_0(x))\nabla u(x)+d_0(x)\nabla \bar u(x)\}\!\cdot\! \{(1-d_0(x))\nabla u(x)+d_0(x)\nabla \bar u(x)\}
$$
$$
+2e^{x_1}a(x)\big(\nabla d_0(x)\big)(\bar u(x)-u(x))\cdot \{(1-d_0(x))\nabla u(x)+d_0(x)\nabla \bar u(x)\}
$$
$$
+e^{x_1}a(x)\big(\nabla d_0(x)\big)(\bar u(x)-u(x))\cdot\big(\nabla d_0(x)\big)(\bar u(x)-u(x))
$$
Using the convexity of the scalar product, this is
$$
\le (1-d_0(x))e^{x_1}a(x)\nabla u(x)\cdot\nabla u(x)+
d_0(x)e^{x_1}a(x)\nabla\bar u(x)\cdot\nabla\bar u(x)
$$
$$
+2e^{x_1}\kappa^{-1}|\nabla d_0(x)| \, |\bar u(x)-u(x)|\, (|\nabla u(x)|+|\nabla \bar u(x)|)
$$
$$
+e^{x_1}a(x)\big(\nabla d_0(x)\big)(\bar u(x)-u(x))\cdot\big(\nabla d_0(x)\big)(\bar u(x)-u(x))\, .
$$
For $x\in\Pi_1$, we have
$$
e^{x_1}a(x)\nabla v_0(x)\cdot\nabla v_0(x)=e^{x_1}a(x)\nabla \bar u(x)\cdot\nabla\bar u(x).
$$
After integrating the former inequality over $\Pi_0\setminus\Pi_1$ and the latter over $\Pi_1$, and summing up, we get
$$
\EE(u,\Pi_0)\le \EE(v_0,\Pi_0)\le \EE(u,\Pi_0\setminus\Pi_1)+ \EE(\bar u,\Pi_0)
$$
$$
+4\kappa^{-1}e^{-L/d} \int\limits_{\Pi_0\setminus\Pi_1} e^{x_1}|\bar u(x)-u(x)|\, (|\nabla u(x)|+|\nabla \bar u(x)|)dx
$$
$$
+4 \kappa^{-1}e^{-2L/d} \int\limits_{\Pi_0\setminus\Pi_1} e^{x_1}( u(x)-u(x))^2\, dx.
$$
After simple rearrangements this yields, using H\"older's inequality,
$$
\EE(u,\Pi_1) \le \EE(\bar u,\Pi_0)+4 \kappa^{-1}e^{-2L/d} \int\limits_{\Pi_0\setminus\Pi_1} e^{x_1}( u(x)-u(x))^2\, dx
$$
\begin{equation}\label{eq1}
+4\kappa^{-1}e^{-L/d}\Big(\int\limits_{\Pi_0\setminus\Pi_1} e^{x_1}( u(x)-u(x))^2\,dx\Big)^{1/2} \Big(\int\limits_{\Pi_0\setminus\Pi_1} e^{x_1}(|\nabla\bar u(x)|^2+|\nabla u(x)|^2)\,dx\Big)^{1/2}
\end{equation}
Our next aim is to estimate the integrals on the right-hand side of the last inequality in terms of the energies of $u$ and $\bar u$.

First, we estimate
$$
\int\limits_{\Pi_0\setminus\Pi_1} e^{x_1}((u(x)-u(x))^2\,dx\le \int\limits_{\Pi_0\setminus\Pi_1} e^{x_1}( u(x)-1)^2\,dx+\int\limits_{\Pi_0\setminus\Pi_1} e^{x_1}(1-u(x))^2\,dx.
$$
Since $(u-1)=0$ for $x_1=L$, we have
$$
(u(x)-1)^2=\Big(\int\limits_{x_1}^L\frac{\partial u}{\partial x_1}(y,z)dy\Big)^2\le 2L\int\limits_{x_1}^L|\nabla u(y,z)|^2dy\, ,
$$
where we used Jensen's inequality. Multiplying this bound by $e^{x_1}$ and integrating over the set $\Pi_0\setminus \Pi_1$, we obtain
\begin{eqnarray}\label{eq2}
\int\limits_{\Pi_0\setminus\Pi_1} e^{x_1}(1-u(x))^2\,dx
&\le&
\int\limits_{\Pi_0\setminus \Pi_1}e^{x_1} 2L\int\limits_{x_1}^L|\nabla u(y,z)|^2dydx \cr
&\le&
\int\limits_{\Pi_0\setminus \Pi_1} 2L\int\limits_{x_1}^L e^{y}|\nabla u(y,z)|^2dydz\cr
&\le&
4L^2\int\limits_{\Pi_0\setminus \Pi_1}e^{x_1}|\nabla u(x)|^2dx
\end{eqnarray}
Similarly, taking into account that $\bar u\not=1$ only for $x_1\le-L+1$, we obtain
$$
\int\limits_{\Pi_0\setminus\Pi_1} e^{x_1}(1-\bar u(x))^2\,dx\le
\int\limits_{\Pi_0\setminus \Pi_1 \cap \{ x_1 \le -L+1\} }e^{x_1}dx\le (e-1)2^{d-1}e^{-L}e^{(d-1)L/d}
= a_d e^{-\frac{L}{d}},
$$
where $a_d =(e-1)2^{d-1}$. Combining the latter bound with (\ref{eq2}) yields
\begin{eqnarray}\label{extra1}
\int\limits_{\Pi_0\setminus\Pi_1} e^{x_1}((u(x)-u(x))^2\,dx
&\le&
a_d e^{-\frac{L}{d}}+4L^2\int\limits_{\Pi_0\setminus \Pi_1}e^{x_1}|\nabla u(x)|^2dx\cr
&\le&
a_d e^{-\frac{L}{d}}+4L^2\kappa^{-1}\int\limits_{\Pi_0\setminus \Pi_1}e^{x_1}a(x)\nabla u(x)\cdot\nabla u(x)dx\cr
&=&
a_d e^{-\frac{L}{d}}+4L^2\kappa^{-1}E(u,\Pi_0 \setminus \Pi_1)
\end{eqnarray}
In the same way, using $|\nabla\bar u(x)|= {\bf 1}_{\{x_1 \le -L+1\}}$, for the last integral on the right-hand side of (\ref{eq1}) we have
\begin{equation}\label{extra2}
\int\limits_{\Pi_0\setminus\Pi_1} e^{x_1}(|\nabla\bar u(x)|^2+|\nabla u(x)|^2)\,dx\le a_d e^{-\frac{L}{d}}+
\Lambda^{-1}E(u,\Pi_0\setminus\Pi_1)
\end{equation}
The following bound for $\EE(\bar u,\Pi_0)$ is straightforward:
$$
\EE(\bar u,\Pi_0)\le \kappa^{-1}\int\limits_{\Pi_0}e^{x_1}|\nabla\bar u(x)|^2 dx \le \kappa^{-1}a_d e^{-\frac{L}{d}}.
$$
>From (\ref{eq1}), (\ref{extra1}), (\ref{extra2}) and the last bound, we derive
the inequality
\begin{eqnarray}\label{eq3}
\EE(u,\Pi_1)
&\le&
 \kappa^{-1}a_d e^{-\frac{L}{d}}\!+\! \kappa^{-1}e^{-2L/d}\big(a_de^{-\frac{L}{d}}
+ 4L^2\kappa^{-1} \EE(u,\Pi_0\setminus \Pi_1)\big)\cr
&\hbox{ }&
+4\kappa^{-1}e^{-L/d}\Big(a_d e^{-\frac{L}{d}}+4L^2 \kappa^{-1}\EE(u,\Pi_0\setminus\Pi_1)\Big)^{1/2}\cr
&\hbox{ }&
\cdot \Big(a_d e^{-\frac{L}{d}}+ \kappa^{-1}\EE(u,\Pi_0\setminus\Pi_1)\Big)^{1/2}\cr
&\le&
 C_{\kappa, d}\big(e^{-\frac{L}{d}}+L^2e^{-L/d} \EE(u,\Pi_0\setminus\Pi_1)\big)\le C_{\kappa, d}\big(e^{-\frac{L}{d}}
+L^2e^{-L/d} \EE(u,\Pi_0)\big)
\end{eqnarray}
Let us now estimate the energy $\EE(u,\Pi_0)$. To this end we denote $G(r,x)=x+[-r,r]^d$, and notice that,
by the standard elliptic estimates (see \cite{kn:GT}),
$$
\|\nabla u\|_{L^2(G(1,x)\cap\Pi)}\le C_1(\kappa, d) \| u\|_{L^2(G(2,x)\cap\Pi)}\le C_2(\kappa, d)
$$
with a constant $C_2(\kappa, d)$ which depends only on $\kappa$ and the dimension.
This implies the bound
$$
\EE(u,\Pi_0)\le C_4(\kappa,d)e^L e^{(d-1)L/d}= C_4(\kappa, d)e^{(2d-1)L/d}.
$$
Together with (\ref{eq3}) this gives
\begin{equation}\label{eq4}
\EE(u,\Pi_1)\le C_5(\kappa,d)L^2e^{(2d-2)L/d}.
\end{equation}
In exactly the same way as (\ref{eq3}), we obtain (for $L^2 > C_5(\kappa, d)$)
$$
\EE(u,\Pi_2)\le C_5(\kappa,d)\big(L^2e^{-L/d}E(u,\Pi_1)+e^{-\frac{L}{d}}\big)\le
C_5(\kappa,d)L^4e^{(2d-3)L/d}.
$$
Iterating this procedure $2d$ times we finally get
$$
\EE(u,\Pi_{2d})\le C_5(\kappa,d)\big(L^{4d - 2}e^{-L/d}+e^{-\frac{L}{d}}\big) \le C_5(\kappa,d)L^{4d}e^{-L/d}.
$$
If $L$ is sufficiently large, then this estimate implies the bound
\begin {equation}\label{LBcontra}
\EE(u,\Pi_{2d})\le e^{-L/(2d)}.
\end{equation}
Clearly, for all sufficiently large $L$ we have $[-L,L]\times[-1,+1]^{d-1}\subseteq \Pi_{2d}$,
and (\ref{LBcontra}) contradicts the lower bound (\ref{lowereng}). We conclude that there is $L_0$ such that for all $L\ge L_0$, we have
$$
u(0,z)\ge \frac{2}{3},\qquad \forall z\in\mathbb R^{d-1}.
$$
Our arguments also ensure that the constant $L_0$ depends only on $\kappa$ and $d$. This completes the proof.
\qed

\begin{lemma} \label{lm:ext} For any $p\geq1$
there exists a constant $C_0$ that depends only on $p$, the dimension, the ellipticity constant $\kappa$
and the $L^\infty$ bound on $V$, such that for all $\l\leq 1$ and $t\geq 1/\l^2$ and for any environment $\o$,
\beqn \label{eq:ext}
E_0^{\l,\,\o}[\,\max\limits_{0\le s\le t}\vert X(s)\vert^p\,]\leq C_0 \lambda^p t^p\,.
\eeqn
\end{lemma}

\medskip\noindent
{\it Proof}

First assume that $V=0$. Letting $\widetilde{X}(t) = \l X\left( \frac{t}{\l^2}\right), t \geq 0$,
we reduce (\ref{eq:ext}) to the following inequality: for all $t\ge 1$ and for all $\l\le 1$,
\beqn \label{eq:ext-scale}
E_0^{\l,\,\o}[\,\max\limits_{0\le s\le t}\vert \widetilde X(s)\vert^p\,]\leq C_0 t^p\,.
\eeqn
Denote $\widetilde T_r=\inf\{s>0\,:\,|\widetilde X(s)|=r\}$.

We now rely on Aronson's lower bound for the Green function ${\cal G}(x,y,t)$  of
the parabolic problem
$$
2\partial_t U-\mathrm{div}(a(x)\nabla U)-a(x)e_1\cdot \nabla U =0\,, \quad U\big|_{|x|=1} =0.
$$
According to \cite{kn:aronson}, Theorems 8 and 9,
for $t\leq 1$ and $x$ with $|x|\le 1/2$ we have
$$
{\cal G}(x,y,t)\geq kt^{-d/2}\exp(-K|x-y|^2/t)
$$
with constants $k>0$ and $K>0$ which only depend on $\kappa$ and $d$.
This implies the bound
$$
P_0^{\l,\,\o}[\widetilde T_1\geq 1]\geq \delta_0>0
$$
which, in turn, yields
\begin{equation}\label{exit_bound}
E_0^{\l,\,\o}\big[e^{-\widetilde T_1}\big]\le 1-\eps_0
\end{equation}
for some $\eps_0>0$ which only depends on $\kappa$ and $d$.

We have
$$
E_0^{\l,\,\o}[\,\max\limits_{0\le s\le t}\vert \widetilde X(s)\vert^p\,]= \int_0^\infty pr^{p-1}P_0^{\l,\,\o}\big[\widetilde T_r\leq t\big]dr =
$$
$$
=\int_0^\infty pr^{p-1}P_0^{\l,\,\o}\big[e^{-\widetilde T_r}\geq e^{-t}\big]dr \leq
e^t\int_0^\infty pr^{p-1}E_0^{\l,\,\o}\big[e^{-\widetilde T_r}\big]dr
$$
Considering the inequality
$$
\widetilde T_r\geq \widetilde T_1+\widetilde T_1\circ\theta_{\widetilde T_1}+
\widetilde T_1\circ\theta_{\widetilde T_2}+\dots+\widetilde T_1\circ\theta_{\widetilde T_{\lfloor r-1\rfloor}},
$$
by the Markov property and (\ref{exit_bound}), for all $t\le 1$ we have
$$
E_0^{\l,\,\o}[\,\max\limits_{0\le s\le t}\vert \widetilde X(s)\vert^p\,]\le
e^t\int_0^\infty pr^{p-1}E_0^{\l,\,\o}\big[e^{-\widetilde T_r}\big] dr\le
e\int_0^\infty pr^{p-1}(1-\eps_0)^{\lfloor r\rfloor}dr\le C(p,\kappa,d).
$$
Using Jensen's inequality and the Markov property again, for all $t\ge 1$ we obtain
$$
E_0^{\l,\,\o}[\,\max\limits_{0\le s\le t}\vert \widetilde X(s)\vert^p\,]\leq
$$
$$
\leq
E_0^{\l,\,\o}\big[\big(\max\limits_{0\le s\le 1}|\widetilde X(s)|+ \max\limits_{1\le s\le 2}|\widetilde X(s)-\widetilde X(1)|+
\ldots+
\max \limits_{\lfloor t\rfloor\le s\le t}|\widetilde X(s)-\widetilde X(\lfloor t\rfloor)|\big)^p\big]
$$
\beqnn
\le (t+1)^{p-1}\Big(E_0^{\l,\,\o}
\big[\max\limits_{0\le s\le 1}|\widetilde X(s)|^p\big] &&+
E_0^{\l,\,\o}\big[\max \limits_{1\le s\le 2}|\widetilde X(s)-\widetilde X(1)|^p\big] \\
&+&\ldots+ E_0^{\l,\,\o}\big[\max\limits_{\lfloor t\rfloor\le s\le t}|\widetilde X(s)-\widetilde X(\lfloor t\rfloor)|^p\big]\Big)
\eeqnn
$$
\le C(p,\kappa,d)(t+1)^p\le C_1(p,\kappa,d)t^p.
$$
Recalling the definition of $\widetilde{X}$, we see that this is equivalent to (\ref{eq:ext}). Hence the proof is complete in the case $V=0$.

To extend the statement  to the case $V\not=0$, we use the time change arguments
from the last part of the proof of Proposition \ref{lm:LR}. We observe that
as in the $V=0$ case, the process $Y^{\l,\omega}$ satisfies estimate (\ref{eq:ext}), and due to (\ref{eq:A}), a similar bound also holds for the process $X^{\l,\omega}$. This completes the proof of the Lemma.
\qed

Lemma \ref{lm:ext} implies the following bound on the effective drift:
\beqn \label{eq:boundondrift}
\vert\ell(\lambda)\vert\leq C_0 \l\,.
\eeqn


\vskip 1cm
\section{Renewal structure}
\setcounter{equation}{0}
\label{sec:V}

As already mentioned in the introduction, $X$ satisfies a Law of Large Numbers
under $\P_0^\l$. The next Proposition is a quantitative version of this convergence.

\begin{proposition} \label{prop:lln}
\beqn\label{eq:lln}
\lim_{\alpha\rightarrow +\infty}
\limsup_{\l\rightarrow 0\,;\, t\rightarrow +\infty \,;\, \l^2t=\alpha}
\left\vert\E_0^\l\left[\frac {X(t)}{\l t}\right] -\frac{\ell(\l)}\l\right\vert=0\,.
\eeqn
\end{proposition}

The proof of Proposition \ref{prop:lln} is based on a renewal argument.
More precisely, (\ref{eq:lln}) follows if we can construct a renewal structure
such that the interval
between two successive renewal times is of order $1/\l^2$.


Our definition of regeneration times is a variant of that in \cite{kn:LS}
where the construction depends on $\l$, whereas in
\cite{kn:LS} $\l$ was fixed.

We shall also heavily rely on the PDE estimates proved in the previous section of the paper.

The first issue we have to address is to check that the approach developed in \cite{kn:LS} applies to our model and does yield the Law of Large Numbers.
Besides in the proof of Proposition \ref{prop:lln} we need sharp estimates on the regeneration times.

\subsection{Construction of  regeneration times}

We recall that $\l$ is chosen small enough. In particular, we assume that $\frac{1}{\l}$ is much larger than the
range $R$ in Assumption 4. We shall also need a constant $l\geq 1$ chosen so that
$c e^{-\kappa_1 l}\leq\frac 12$ where $c$ and $\kappa_1$ are the constants appearing in
Lemma \ref{lm:leftexit}.

We set $R(\l) : = \frac{l}{\l}$. We now follow the construction of \cite{kn:LS}, replacing $R$ in his construction with $R(\l)$.
For details of this construction, proofs (and for pictures), we refer to \cite{kn:LS}.
We first have to enlarge the probability space by adding an auxiliary sequence $(Y_k)_{k\geq 0}$ of i.i.d. Bernoulli random variables.
Denote $B_R(x)$ the ball with center $x$ and radius $R$.
Let $U^x: = B_{6R(\l)}(x + 5R(\l)e_1)$, $B^x: = B_{R(\l)}(x+ 9R(\l)e_1)$, and
let
\beq\label{exittimedef}
T_{\rm{exit}, U^x} : = \inf\{s\geq 0: X(s) \notin U^x \}
\eeq
be the exit time from $U^x$.
We consider the corresponding transition density $p_{\l, \o, U^x}(s,x,y)$ which is defined by
$P^{\l, \o}_x[X(s) \in G, T_{\rm{exit}, U^x} > s] = \int\limits_G p_{\l, \o, U^x}(s,x,y) dy$, for all open sets $G\subseteq U^x$.
We will need the following bound for this transition density.

\begin{lemma} \label{lm:btrans}
There is  some $\delta > 0$ (depending on $V$, $\sigma$ and $d$) such that
\beq\label{btrans}
p_{\l, \o, U^x}(1/\l^2,x,y)
\geq \frac{2\delta}{|B_{R(\l)}|}, \hbox { for all } x\in \R^d, \,  y \in B^x\,
\hbox{ and } \l \leq 1 .
\eeq
\end{lemma}
{\it Proof}

Again we begin with the case $V=0$.
After rescaling $t/\l^2\to t$, $(x/\l,\, y/\l)\to(x,\,y)$ the required bound is an immediate consequence of Theorems 8,\,9 in \cite{kn:aronson}.

If $V\not=0$, then the desired lower bound is an immediate consequence of the following statement:

\begin{lemma}\label{l_harhar}
Let a function $\rho(x)$ and a symmetric matrix $\{\alpha_{ij}(x)\}$ satisfy the estimates
$$
\varkappa\leq\rho\le \varkappa^{-1},\qquad \varkappa I\leq\alpha\le \varkappa^{-1}I,\qquad\varkappa>0,
$$
and denote by $G_B(t,x,y)$ a solution to the following parabolic problem
\begin{equation}\label{aux_para_e}
\begin{array}{rl}
\displaystyle
\rho(x)\partial_t G_B(t,x,y)-\mathrm{div}_x\!\!&\!\!\big(\alpha(x)\nabla_x G_B(t,x,y)\big)=0,\\[3mm]
G_B\big|_{x\in\partial B}=0,\ \ &\ \ G_B(0,x,y)=\delta(x-y)
\end{array}
\end{equation}
with $B$ being an open ball centered at the origin, and $y\in B$. Then for any ball $B_0$ such that
$\overline B_0\subset B$, the following inequality holds
\begin{equation}\label{bnd_for_l52}
G_B(1,x,0)\ge C\qquad\hbox{for } x\in B_0,
\end{equation}
with a constant $C$ which only depends on $\varkappa$, $d$, $B$ and $B_0$.
\end{lemma}

As in the case $V=0$ the estimate (\ref{btrans}) can be obtained from (\ref{bnd_for_l52}) by scaling $t/\l^2\to t$, $(x/\l,\, y/\l)\to(x,\,y)$.

\medskip
\noindent
{\it Proof of Lemma \ref{l_harhar}}.\ \
Enlarging if necessary the ball $B_0$ we assume without loss of generality, that $B_0$ contains the origin.

Consider an auxiliary spectral problem in $B$
$$
-\mathrm{div}\big(\alpha(x)\nabla \Psi(x)\big)=\rho(x)\nu\Psi(x), \qquad\Psi\big|_{\partial B}=0.
$$
By means of the minimax principle one can check that the principal eigenvalue $\nu_1$  satisfies the estimate $0<\nu_1<C_1$. The principal eigenfunction $\Psi_1$ is known to be positive in $B$. Assuming the normalization
$$
\int_B\Psi_1(x)dx=1,
$$
by the Harnack inequality and H\"older continuity arguments (see \cite{kn:GT}) we conclude that
\begin{equation}\label{haba}
\Psi_1(x)\leq C_2\quad\hbox{in }B,\qquad \Psi_1(x)\geq C_3\quad\hbox{in }B_0,
\end{equation}
where the constant $C_1$ and $C_2$ depend only on $\varkappa$, $d$ and $B$, and
$C_3$ also depends on $B_0$.
Clearly, the function $e^{-\nu_1t}\Psi_1(x)$ solves problem (\ref{aux_para_e})
with the initial condition $\Psi_1$. Therefore,
$$
e^{-\nu_1/2}\Psi_1(x)=\int_B G_B(1/2,x,y)\Psi_1(y)dy.
$$
Making use of (\ref{haba}) and the upper bound for $\nu_1$ we derive the inequality
$$
\int_B G_B(1/2,0,y)dy\ge C_4.
$$
Considering the symmetry  of the operator with respect to the weighted measure $\rho\,dx$, we have $\rho(x)G_B(t,x,y)=\rho(y)G_B(t,y,x)$.
It readily follows from the results of \cite{kn:GuCr} that the function $G(1/2,y,0)$ satisfies the upper bound $G_B(1/2,y,0)\le C_5$ in $B$ with a constant $C_5$ which only depends
on $\kappa$ and $d$.
Consequently, there is a smaller ball $B_1$ centered at the origin, $\overline B_1\subset B$, such that
$$
\int_{B_1} G_B(1/2,y,0)dy\ge C_4/2,
$$
the radius of $B_1$ depends only on $B$, $\kappa$ and $d$. This yields
\begin{equation}\label{intbou}
\sup\limits_{y\in B_1}G_B(1/2,y,0)\ge C_6.
\end{equation}
Without loss of generality we suppose that $B_0\subset B_1$.

According to \cite[Theorem 1-5]{kn:GuCr}, the function $G_B(1,y,0)$ satisfies the following version of the Harnack inequality in $B_1$:
$$
\inf\limits_{y\in B_1}G_B(1,y,0)\ge C_7(\varkappa,d)\sup\limits_{y\in B_1}G_B(1/2,y,0).
$$
This estimate combined
with (\ref{intbou}) yields (\ref{bnd_for_l52}). This completes the proof of Lemmata \ref{l_harhar} and \ref{lm:btrans}.
\qed

\bigskip
Due to (\ref{btrans}), we can give the following coupling construction.
Let $(\FF_t)_{t \geq 0}$ be the filtration generated by $(X(t))_{t \geq 0}$ and
${\mathcal S}_m : = \sigma(Y_0, \ldots , Y_m)$.

We denote $\theta_m^\l$ the rescaled shift operator defined by
$$
\theta_m^\l((X(s)_{s\geq 0}) = (X({\l^{-2}m+s}))_{s\geq 0}.
$$
These shift operators $\theta^\l_m$, $m \in \N$, are extended in the obvious way:
$\theta^\l_m((X(s))_{s\geq 0},(Y_k)_{ k\geq 0}) = ((X(\l^{-2}m+s))_{s\geq 0},$ $(Y_{m+k})_{ k\geq 0})$.

\begin{proposition} \label{prop:coupling}
There exists, for every $\l$, $\o$ and $x$, a probability measure $\widehat{P}^{\l, \o}_x $ on the enlarged probability space
such that, with $\delta$ from (\ref{btrans}), \\
(i) The law of $(X(t))_{t \geq 0}$
under $\widehat P^{\l, \o}_x $ is $P^{\l, \o}_x $, and the sequence $(Y_k)_{k \geq 0}$ is
a sequence of i.i.d. Bernoulli variables with success probability $\delta$ under $\widehat{P}^{\l, \o}_x $.\\
(ii) Under $\widehat P^{\l, \o}_x $,
$(Y_n)_{n \geq m}$
is independent of $\FF_{\l^{-2}m} \times {\mathcal S}_{m-1}$,
and conditioned on
$\FF_{\l^{-2}m} \times {\mathcal S}_m$, $X\circ \theta^\l_m$ has the same law as $X$ under
$\widehat{P}^{\l, \o}_{X(\l^{-2}m),Y_m} $, where
$\widehat{P}^{\l, \o}_{x,y} $ denotes the conditioned law $\widehat{P}^{\l, \o}_x[\cdot | Y_0 =y]$, (for $y \in \{0,1\}$).
\\
(iii) $\widehat{P}^{\l, \o}_{x,1} $-almost surely, $X(t) \in U^x$ for $t \in [0,\l^{-2}]$ and
the distribution of $X(\l^{-2})$ under $\widehat{P}^{\l, \o}_{x,1} $ is the uniform distribution on $B^x$.
\end{proposition}
We refer to \cite{kn:LS} for the proof.

We will now introduce random times $N_k\in \l^{-2}\mathbb Z^+$ for which $Y_{\l^2N_k}=1$ and for which the process
$(e_1\cdot  X(t))_{t \geq 0}$ essentially reaches a local maximum (within a variation of $R(\l)$). The first regeneration time $\tau_1$ will be the first time $N_k +\l^{-2}, k\geq 1$ such that $(e_1\cdot X(t))_{t \geq 0}$ never goes below
$e_1\cdot X(N_k+\l^{-2})-R(\l)$ after $N_k +\l^{-2}$.
In order to define $N_k$, we will first consider  stopping times $\widetilde N_k \in \l^{-2}\mathbb Z^+$ which are essentially the times when $(e_1\cdot X(t))_{t \geq 0}$ reaches local maxima (also within a variation of $R(\l)$). Then, $N_1$ will be the first $\widetilde N_k$ with $Y_{\l^2\widetilde N_k} =1$.

Let
\beq\label{Mdef}
M(t): = \sup\{e_1\cdot( X(s)-  X(0)): 0 \leq s \leq t\}
\eeq
For $a> 0$, define the stopping times $V_k^\l(a), k \geq 1$, as follows. Recall that $T_L=\inf\{t\,: \, \tr{e_1}\,X(t)=L\}$, and define
\beq
V_0^\l(a) : =  T_{\l^{-1}a}, \quad V_{k+1}^\l(a) : = T_{M(\lceil V_k^\l(a)\rceil_\l)+ R(\l)}, \quad k \geq 1;
\eeq
here and later on $\lceil r\rceil_\l$ stands for the $\min\{n\in\mathbb \l^{-2}\mathbb Z\,:\, r\le n\}$.
Then
\beq
\widetilde N_1^\l(a): =\inf\left\{\big\lceil V_k^\l(a)\big\rceil_\l: k\geq 0, \sup\limits_{s \in [V_k^\l(a), \lceil V_k^\l(a) \rceil_\l ]}\left|e_1\! \cdot \Big( X(s) - X\big(V^\l_k(a))\Big)\right| \leq \frac{ R(\l)}{2}\right\},
\eeq
\beq
\widetilde N^\l_{k+1}(a): = \widetilde N_1^\l(3\l R(\l))\circ \theta^\l_{\l^2\widetilde N_k^\l(a)} + \widetilde N_k^\l(a), \quad k\geq 1\, ,
\eeq
\beq
N_1^\l(a) : = \inf\left\{\widetilde N_k^\l(a): k \geq 1, Y_{\l^2\widetilde N_k^\l(a)} =1\right\}\, ,
\eeq
(we will see later that $\widetilde N_k^\l(a) < \infty$, for all $k$). The random times $\l^2\widetilde N_k^\l(a)$ are integer-valued and $\sup\limits_{s\leq \widetilde N_k^\l(a)}e_1\cdot (X(s) - X(\widetilde N_k^\l(a))) \leq R(\l)$.
We next define random times $S_1, J_1$ and $R_1$ as follows.
\beq
S_1^\l: = N_1^\l(3\l R(\l)) +\l^{-2}, \quad J_1^\l : = S_1^\l+ T_{- R(\l)} \circ \theta^\l_{\l^2 S_1^\l}, \quad R_1^\l : = \lceil J_1^\l\rceil_\l = S_1^\l + D\circ \theta^\l_{\l^2S_1^\l} ,
\eeq
where
\beq\label{Ddef}
D: = \lceil T_{-R(\l)}\rceil_\l\,  .
\eeq
Now we proceed recursively:
\beq
N_{k+1}^\l = R_k^\l + N_1^\l(a_k)\circ \theta^\l_{\l^2R_k^\l} \quad \hbox{ with } a_k = \l\big(M(R_k^\l) - e_1 \cdot(X(R_k^\l)- X(0)) +R(\l)\big)
\eeq
and
$$
S_{k+1}^\l: = N_{k+1}^\l +\l^{-2}, \quad J^\l_{k+1} : = S^\l_{k+1}+ T_{R(\l)} \circ\theta^\l_{\l^2S^\l_{k+1}} ,
\quad R_{k+1}^\l : = \lceil J_{k+1}^\l\rceil_\l = S_{k+1}^\l + D\circ \theta^\l_{\l^2S^\l_{k+1}} \, .
$$
Note that for all $k$, the $\FF_t \times {\mathcal S}_{\l^2\lceil t\rceil_\l}$- stopping times $\l^2N^\l_k,\ \l^2S^\l_k$ and $\l^2R^\l_k$ are integer-valued (the value $+\infty$ is possible). By definition, we have $\l^{-2} \leq N^\l_1 \leq S^\l_1 \leq J^\l_1\leq R^\l_1 \leq N^\l_2 \leq S^\l_2 \leq J^\l_2 \leq R^\l_2 \leq N^\l_3 \ldots \leq \infty$.
The first regeneration time $\tau_1$ is defined as
\beq
\tau_1 : = \inf\{S_k^\l: S_k^\l < \infty, \, \, R_k^\l = \infty \} \leq \infty\, .
\eeq
By definition, $\l^2\tau_1$ is integer-valued and $\tau_1 \geq 2\l^{-2}$ (since $N_1^\l \geq \l^{-2}$). We see that on the event $\tau_1 < \infty$ it holds
$$
e_1\cdot X(s) \leq e_1\cdot X(\tau_1-\l^{-2})+ R(\l)\leq  e_1\cdot X(\tau_1) - 7R(\l),\quad \hbox{for }s \leq \tau_1 -\l^{-2}, \ \ \widehat P^{\l,\o}_x-\hbox{a.s.},
$$
see also Proposition \ref{prop:coupling}, \ i.e. $(X(s))_{s\leq \tau_1 -\l^{-2}}$ remains in the halfspace
$\{z \in \R^d: e_1\cdot z \leq  e_1\cdot X(\tau_1) -7 R(\l)\}$. On the other hand, since the process $(e_1\cdot X(t))_{t \geq 0}$ never goes below $e_1 \cdot X(\tau_1) - R(\l)$ after $\tau_1$, $\widehat P^{\l,\o}_x$-a.s., $( X(t))_{t > \tau_1}$ remains in the halfspace
$\{z\in \R^d: e_1 \cdot z \geq e_1\cdot X(\tau_1) - R(\l)\}$.

In \cite{kn:LS}, it is proved that $\tau_1 < \infty$ if and only if the process is transient in direction $e_1$.
More precisely, define the annealed law

\beq\label{eq:annhat}
\widehat\P^{\l}_x[A]: =\int d\Q(\o)\int d\widehat P^{\l,\o}_x(w)\1_A(\o,w)\, .
\eeq

\begin{proposition}\label{finitreg}
The following statements are equivalent:\\
(i) $\tau_1 < \infty$, $\widehat\P^{\l}_0$-a.s.\\
(ii) $e_1\cdot X(t) \to \infty$, $\P^{\l}_0$-a.s. \\
(iii) $\P^{\l}_0[ D= \infty] > 0$.
\end{proposition}
For the proof, we refer to \cite{kn:LS}, Lemma 2.3 and Proposition 2.7.
Later on we will need stronger results than  those of Proposition \ref{finitreg}
involving sharp bounds on $\tau_1$ as $\l\to0$.

We set $\tau_0=0$ for convenience.
The next theorem (Theorem 2.5 in \cite{kn:LS}) gives the renewal structure which is crucial to establish
(for fixed values of $\l$) the Law of Large Numbers and the Central Limit Theorem.

\begin{theorem}\label{thm:renew}
Assume that $\tau_1 < \infty$, $\widehat\P^{\l}_0$-a.s.
Then, under the measure $\widehat\P^{\l}_0$, the random variables $Z_k: = \left((X((\tau_k +t)\wedge(\tau_{k+1}-\l^{-2}))- X(\tau_k))_{t \geq 0}, X(\tau_{k +1})- X(\tau_k), \tau_{k+1} - \tau_k \right), k \geq 0$ are independent. Furthermore, the random variables $Z_k, k\geq 1$, are i.i.d.
under $\widehat\P^{\l}_0$ and have the same law as $Z_0$ under $\widehat\P^{\l}_0[\,  \cdot\, | D=\infty]$.
\end{theorem}

Note that the renewal structure is proved for the trajectory between the times $\tau_k$ and $\tau_{k+1}-\l^{-2}$, but we have a good
control over the trajectory between the times $\tau_{k+1}-\l^{-2}$ and $\tau_{k+1}$  : since
$Y_{\l^2\tau_{k+1} -1} = 1$, then $X(s)\in  U^{X_{\tau_{k+1} -1/\l^{2}}}$, for all $s \in [\tau_{k+1} -\l^{-2}, \tau_{k+1}]$, i.e. the trajectory remains in a ball of radius $6R(\l)$.

Let
\beq
K=\inf\{k\geq 1\,:\, S_k^\l<\infty\hbox{\, and\,} R_k^\l=\infty\}.
\eeq
Then $\tau_1 = S^\l_K$.
The points $X(S_1^\l),\, X(S_2^\l),\, \ldots $ are {\it ladder points} of the process.
The idea of such a decomposition
of the path goes back to \cite{kn:kesten} and was first turned effective for multi-dimensional
random walks in random environments in \cite{kn:SZ}.

\medskip

\begin{lemma} \label{lm:tau1isfinite}
The following statements hold:

\medskip\noindent
(i) $S_1^\l$ is $\widehat \P_0^\l$-almost surely finite and, for all $k$, $S_k^\l$ is $\widehat \P_0^\l$-almost surely finite
on the event $R^\l_{k-1}<\infty$. \\
(ii) $\tau_k$ is $\widehat \P_0^\l$ almost surely finite for all $k \geq 1$.\\
(iii) $\widehat \E_0^\l[\tau_k] < \infty$ for all $k \geq 1$.
\end{lemma}

\noindent
{\it Proof}
Part (i). Due to Lemma \ref{lm:rightexit}, $V^\l_k(a)< \infty$, $\widehat \P_x^\l$-a.s for all $k$ and all $x$.
We have (as in Lemma \ref{lm:btrans} this is a consequence of the Aronson-Nash lower bounds
for the Green function of a parabolic equation)
\beq\label{unifbdd}
\inf\limits_{\l \leq 1}\inf\limits_{ x \in \R^d} \inf\limits_\o \widehat P_x^{\l, \o}\left[\max\limits_{ s \leq 1/\l^2}|e_1\cdot X(s)|\leq \frac{ R(\l)}{2}\right] = \widetilde\delta > 0\, ,
\eeq
and this implies that $\widetilde N^\l_k(a) < \infty$ for all $a,\,k$, $\widehat \P_x^\l$-a.s for all $x$. Due to Proposition \ref{prop:coupling}, we conclude that $N_k^\l(a) < \infty$ for all $a,\,k$, $\widehat \P_x^\l$-a.s for all $x$.

\medskip\noindent
Part (ii). It follows from Lemma \ref{lm:leftexit}
that under proper choice of $l$ in the definition of ${R(\l)}$, the following bound holds
\begin{equation}\label{K=1bound}
\widehat \P_0^\l[K=1]=\widehat \P_0^\l[R^\l_1=\infty]\geq 1-ce^{-\kappa_1 l}\geq\frac 12\, .
\end{equation}
Lemma \ref{lm:leftexit} and the Markov property
applied at time $S_{k+1}^\l$ also imply that
$$\widehat \P_0^\l[R_{k+1}^\l=\infty \vert R_k^\l<\infty]\geq \frac 12\,.$$

Thus we get that
\begin{equation}\label{Kbound}
\widehat \P_0^\l[K\geq k] \leq 2^{-k+1}
\end{equation}
and
$$\widehat \P_0^\l[K=\infty]=0\,.$$
Together with part (i), this implies part (ii).

\medskip\noindent
We now turn to part (iii). The next lemma gives a bound for the tail of the random variable $\tau_1$, which will be sufficient to
guarantee that $\tau_1$ has finite expectation under $\widehat \P_0^\l$.

\begin{lemma}\label{l:tailbound}
There exist constants $c_3 > 0$ and $c_4 > 0$ such that
\begin{equation}\label{tailtaufirst}
\hbox{for all  $\l \leq 1$ and  $t > 0$,} \quad \widehat \P_0^{\l}[\, \tau_1\geq \l^{-2} t\, ] \leq c_3 e^{-c_4t}\, .
\end{equation}
The same tail bound holds for the differences $\tau_{k+1} - \tau_k$ for all $k\ge 1$:
\begin{equation}\label{tailtau}
\hbox{for all $\l \leq 1$ and  $t > 0$,} \quad
\widehat \P_0^{\l}[\, \tau_{k+1} -
\tau_k\geq \l^{-2} t]
\leq c_3 e^{-c_4t}\, .
\end{equation}
\end{lemma}

\noindent
{\it Proof of Lemma \ref{l:tailbound}}.
It suffices to show (\ref{tailtaufirst}). Then, (\ref{tailtau}) follows since
$$
\widehat \P_0^\l[\, \tau_{k+1} - \tau_k\geq \l^{-2} t] = \widehat \P_0^\l[\, \tau_1\geq \l^{-2} t| D= \infty]
$$
(see Theorem \ref{thm:renew}) and $\widehat \P_0^\l\left[ D = \infty \right]\geq \frac 12$, see
(\ref{K=1bound}). To show \eqref{tailtaufirst}, we claim that the following stronger statement holds:
There exist constants $c_3 > 0$ and $c_4 > 0$ such that for all $\l \leq 1$, $t > 0$ and all $\omega$,
\begin{equation}\label{tailtaufirstquenched}
\widehat P_0^{\l, \omega}[\, \tau_1\geq \l^{-2} t\, ] \leq c_3 e^{-c_4t}\, .
\end{equation}
Since one can follow the proof of Corollary 4.10 in \cite{kn:LS}, we give only a sketch of the proof of \eqref{tailtaufirstquenched}. \\
{\sl Step 1} (corresponds to Proposition 4.7 in \cite{kn:LS}). Recall \eqref{Mdef}. There exists a constant $c_5 > 0$ such that
\begin{equation}\label{Maxfini}
\sup\limits_{\omega} \widehat E_0^{\l, \omega}\left[\exp(c_5 \l e_1 \cdot M(T_{-R(\l)}))\1(T_{-R(\l)} < \infty)\right] < \infty\, .
\end{equation}
Proof of \eqref{Maxfini}: Due to Lemma \ref{lm:leftexit},
$\widehat P_0^{\l, \o}[\l e_1 \cdot M(T_{-R(\l)}) > t, T_{-R(\l)}< \infty] \leq c e^{-\kappa_1(l+t)}$ and this 
implies \eqref{Maxfini}.\\
{\sl Step 2} Follow the proof of Theorem 4.9 in \cite{kn:LS} to obtain that there is a  constant $c_6 > 0$ such that for all $\l \leq 1$,
\begin{equation}
\sup\limits_{\o} \widehat E_0^{\l, \o}\left[\exp(c_6 \l e_1 \cdot X_{\tau_1})\right] < \infty\, .
\end{equation}
{\sl Step 3}
Take $t > \frac{6l}{\kappa_2}$ and $u = \frac{\kappa_2}{2}t$ where $\kappa_2$ is the constant from Lemma \ref{lm:rightexit}. Then,
\begin{equation}\label{casesschlau}
\widehat P_0^{\l, \o}[\, \tau_1\geq \l^{-2} t\, ] \leq \widehat P_0^{\l, \o}[\, \tau_1\geq \l^{-2} t, e_1\cdot X_{\tau_1}< \l^{-1}u - 3 R(\l) ] + \widehat P_0^{\l, \o}[\,e_1\cdot X_{\tau_1}\geq \l^{-1}u - 3 R(\l)]\, .
\end{equation}
Due to Step 2, the second term on the r.h.s. of \eqref{casesschlau} is $\leq e^ {-c_7 u}$ for some $c_7 > 0$.
Turning to the first term in \eqref{casesschlau}, note that since
$$
\sup\limits_{s \leq \tau_1}e_1 \cdot (X_s - X_{\tau_1}) < 3R(\l),
$$
see the definition of $U^x$ and $B^x$ before (\ref{exittimedef}), $\widehat P_0^{\l, \o}[\, \tau_1\geq \l^{-2} t, e_1\cdot X_{\tau_1}< \l^{-1}u - 3 R(\l) ] \leq P_0^{\l, \o}[\, T_{\l^{-1} u} \geq \l^{-2} t]  \leq Ce^{-\kappa_2 t +u}$, where the last inequality comes from Lemma \ref{lm:rightexit}.\qed

As mentioned before, the regeneration structure implies a law of large numbers for fixed $\l$.
\begin{proposition} \label{prop:estX}
We have, for each $\l > 0$,
\begin{equation}\label{lawoflargenumb}
\lim\limits_{t \to \infty} \frac{X(t)} {t} = \frac{\widehat \E_0^\l \left[X(\tau_2) - X(\tau_1)\right]}{ \widehat \E_0^\l \left[\tau_2 - \tau_1\right]}\quad \P_0^\l \hbox{- a.s.}
\end{equation}
As a consequence,
\begin{equation} \label{renewdiff}
\widehat \E_0^\l \left[X(\tau_2) - X(\tau_1)\right] = \widehat \E_0^\l \left[\tau_2 - \tau_1\right]\ell(\l) \, .
\end{equation}
\end{proposition}

{\it Proof}

Theorem \ref{thm:renew} implies, with the ergodic Theorem, that if $\widehat \E_0^\l \left[\tau_2 - \tau_1\right]$ is finite, then
$\lim\limits_{t \to \infty} \frac{X(t)} {t} $ exists, $\P_0^\l $-a.s. and (\ref{lawoflargenumb}) holds true.
Once again, we refer to \cite{kn:LS} for the details.
\qed

\subsection{Estimates on the regeneration times}
We now show that under $\widehat P^\l_0$, $\tau_1$ and $\tau_2 - \tau_1$ are of order $\l^{-2}$. More precisely,
\begin{lemma}\label{lm:tau}
We have
\beqn\label{eq:tau}
\limsup_{\l\rightarrow 0} \l^4 \widehat \E_0^\l[\,\tau_1^2]<\infty \quad \hbox{\rm and}\quad
\limsup_{\l\rightarrow 0} \l^4 \widehat \E_0^\l[\,(\tau_2 -\tau_1)^2]<\infty\,.
\eeqn
As a consequence,
\beqn\label{eq:*}
\limsup_{\l\rightarrow 0} \l^2 \widehat \E_0^\l[\,\tau_1]<\infty \quad  \hbox{\rm and}\quad
\limsup_{\l\rightarrow 0} \l^2 \widehat \E_0^\l[\,\tau_2 - \tau_1]<\infty\,,\eeqn
and
\beqn\label{eq:**}
\limsup_{\l\rightarrow 0} \frac {\widehat \E_0^\l[\,\tau_1^2]}{\widehat \E_0^\l[\,\tau_1]^2}<\infty \quad \hbox{\rm and}\quad \limsup_{\l\rightarrow 0} \frac {\widehat \E_0^\l[\,(\tau_2 - \tau _1)^2]}{\widehat \E_0^\l[\,\tau_2 -\tau_1]^2}<\infty\,.\eeqn
\end{lemma}

{\it Proof}

Observe that (\ref{eq:*}) directly follows from (\ref{eq:tau}) and H\"older's inequality.

To deduce  (\ref{eq:**}) from (\ref{eq:tau}) it suffices to prove a lower bound on
$\widehat \E_0^\l[\,\tau_1]$ (or $\widehat \E_0^\l[\,\tau_2 - \tau_1]$, respectively) of the order $\l^{-2}$. But note that, since $\tau_1\geq T_{l/\l}$, we have
\beqnn
 \l^2\widehat \E_0^\l[\,\tau_1]
\geq\widehat \P_0^\l[\, \tau_1\geq \l^{-2}]
\geq \P_0^\l[\, T_{l/\l}\geq \l^{-2}]\,.\eeqnn

We next use the Girsanov formula (\ref{eq:gir}):
$$
\P_0^\l[\, T_{l/\l}\geq \l^{-2}]
= \int E[\, \1(T_{l/\l}\geq \l^{-2})
e^{\l\bB (\frac 1{\l^2}) -\frac {\l^2}2 \langle\bB\rangle(\frac 1{\l^2})}]\, d\Q\,.$$

Thus the invariance principle implies that this last quantity has a positive limit: namely
it converges to
$E[\, \1(T_l\geq 1)
e^{Z -\frac 12 E(Z^2)}]\,,$
where $Z$ is some Gaussian random variable and $T_l$ is the hitting time of level $l$ by some
Brownian motion. This last expectation cannot be $0$.
The same argument applies to $\tau_2 - \tau_1$.

Proof of (\ref{eq:tau}): due to (\ref{tailtaufirst}), we have
$\widehat\P_0^\l[\, \tau_1^2 \geq t\l^{-4}] =
\widehat \P_0^\l[\, \tau_1\geq \sqrt{t}\l^{-2}]\leq c_3e^{-c_4 \sqrt{t}}$.
In the same way, due to
(\ref{tailtau}), we have
$\widehat\P_0^\l[\, (\tau_2 -\tau_1)^2 \geq t\l^{-4}] \leq c_3e^{-c_4 \sqrt{t}}$,
and (\ref{eq:tau}) follows. \qed

We can now show the corresponding bounds for the regeneration distances.
\begin{lemma}\label{lm:dist}
We have
\begin{equation}\label{distmomentfirst}
\limsup_{\l\rightarrow 0} \l^2 \widehat \E_0^\l \left[\vert X(\tau_1)\vert^2\right] < \infty
\end{equation}
and
\begin{equation}\label{distmoment}
\limsup_{\l\rightarrow 0} \l^2 \widehat\E_0^\l \left[\vert X(\tau_2) - X(\tau_1)\vert^2\right] < \infty\, .
\end{equation}
\end{lemma}

We note for further reference that, as a consequence of (\ref{distmomentfirst}) and (\ref{distmoment}),
\begin{equation} \label{kdistmoment}
\sup\limits_k\frac{1}{k^2}\limsup_{\l\rightarrow 0} \l^2 \widehat\E_0^\l \left[\vert X(\tau_k)\vert^2\right] < \infty\, .
\end{equation}
Moreover,  note that due to (\ref{eq:boundondrift}),
\begin{equation}\label{ellbdd}
\limsup_{\l\rightarrow 0} \frac{\vert \ell(\l)\vert}{\l} < \infty\, .
\end{equation}
(This also follows from (\ref{distmoment}) together with (\ref{renewdiff}), (\ref{eq:*}) and (\ref{eq:**})).

\medskip\noindent
{\it Proof of Lemma \ref{lm:dist} }
To show (\ref{distmomentfirst}), note that

\begin{eqnarray}
\widehat\E_0^\l \left[X(\tau_1)^2\right] &=&
\sum\limits_{k=0}^\infty \widehat \E_0^\l \left[X(\tau_1)^2 \1(k\leq \l^2 \tau_1 < k+1) \right]\cr
&\leq& \sum\limits_{k=0}^\infty \widehat \E_0^\l \left[\max\limits_{t \leq (k+1)/\l^2} |X(t)|^2 \1(k\leq \l^2\tau_1)\right]\cr
&\leq&  2\sum\limits_{k=0}^\infty \left(\E_0^\l \left[\max\limits_{t \leq (k+1)/\l^2} |X(t)|^4\right]\right)^{1/2}
\Big(\widehat \P_0^\l \left[\l^2\tau_1 \geq k\right]\Big)^{1/2}\cr
&\leq& 2\sum\limits_{k=0}^\infty \left(\l \frac{k+1}{\l^2}\right)^2 \left(c_3  e^{-c_4 k} \right)^{1/2} \leq c_8 \l^{-2}\nonumber
\end{eqnarray}
for some constant $c_8>0$, where we used (\ref{eq:ext}) and (\ref{tailtaufirst}).
Now, (\ref{distmoment}) follows from (\ref{distmomentfirst}), since
$$ \widehat\E_0^\l \left[\vert X(\tau_2) - X(\tau_1)\vert^2\right] = \widehat \E_0^\l \left[\vert X(\tau_1)\vert^2 \vert
D= \infty\right]
$$
(see Theorem \ref{thm:renew}) and $\widehat \P_0^\l\left[ D = \infty \right]\geq \frac 12$, for all $\l$, see
(\ref{K=1bound}). \qed

\subsection{Proof of Proposition \ref{prop:lln}}

Let
$$
n(t) : = \left\lfloor \frac{t}{\widehat \E_0^\l [\tau_2 - \tau_1]}\right\rfloor \, .
$$
We write
\begin{equation}\label{xdecomp}
X(t) = X(\tau_{n(t)}) + (X(t) - X(\tau_{n(t)})\, .
\end{equation}
We will show that
\begin{equation}\label{truelim}
\lim_{\alpha\rightarrow +\infty}\limsup_{\l\rightarrow 0\,;\, t\rightarrow +\infty \,;\, \l^2t=\alpha}
\left\vert \frac{1}{\l t} \widehat\E_0^\l [X(\tau_{n(t)})] - \frac{\ell(\l)}{\l} \right\vert = 0
\end{equation}
and that
\begin{equation}\label{errorterm}
\lim_{\alpha\rightarrow +\infty}\limsup_{\l\rightarrow 0\,;\, t\rightarrow +\infty \,;\, \l^2t=\alpha}
\left\vert \frac{1}{\l t}\widehat\E_0^\l\left[ X(t) - X(\tau_{n(t)})\right]\right\vert =0\, .
\end{equation}
(\ref{truelim}) and (\ref{errorterm}) then imply
\beqn
\lim_{\alpha\rightarrow +\infty}
\limsup_{\l\rightarrow 0\,;\, t\rightarrow +\infty \,;\, \l^2t=\alpha}
\left\vert\E_0^\l\left[\frac {X(t)}{\l t}\right] -\frac{\ell(\l)}\l\right\vert=0\, ,
\eeqn
i.e. Proposition \ref{prop:lln}.

To show (\ref{truelim}), note that (recalling $\tau_0=0$),
\begin{equation}\label{writsums}
\left\vert \frac{1}{\l t} \widehat\E_0^\l [X(\tau_{n(t)})] - \frac{\ell(\l)}{\l} \right\vert
=
\left\vert \frac{1}{\l t} \widehat\E_0^\l \left[\sum\limits_{k=1}^{n(t) } (X(\tau_k) - X(\tau_{k-1})) \right]
 - \frac{\ell(\l)}{\l}\right\vert
\end{equation}
Using formula (\ref{renewdiff}) in Proposition \ref{prop:estX}, we rewrite and estimate the r.h.s. of (\ref{writsums}) as
$$
\left\vert \frac{1}{\l t} \widehat \E_0^\l \left[X(\tau_1)\right] + \frac{1}{\l t} \left(n(t) -1\right)
\widehat \E_0^\l \left[X(\tau_2)- X(\tau_1)\right] - \frac{\ell(\l)}{\l}\right\vert
$$
\begin{equation}\label{firstisdiff}
\leq
\left\vert \frac{1}{\l t}  \widehat\E_0^\l \left[X(\tau_1)\right] \right\vert
+ \left\vert \frac{\ell(\l)}{\l} \left(\frac{\widehat\E_0^\l [\tau_2 - \tau_1]}{t} \left(
\left\lfloor\frac{t}{\widehat\E_0^\l [\tau_2 - \tau_1]}\right\rfloor-1\right) -1\right)\right\vert \, .
\end{equation}

Due to (\ref{distmomentfirst}), the first term in the right hand side of (\ref{firstisdiff})
is of order $\alpha^{-1}$ and therefore satisfies
$$
\lim_{\alpha\rightarrow +\infty}\limsup_{\l\rightarrow 0\,;\, t\rightarrow +\infty \,;\, \l^2t=\alpha}\vert \frac{1}{\l t}  \widehat\E_0^\l \left[X(\tau_1)\right] \vert = 0\, .
$$
Now consider the second term in (\ref{firstisdiff}).
>From (\ref{eq:*}), we know that $\widehat\E_0^\l [\tau_2 -\tau_1]/t$ is of order $\alpha^{-1}$.
Therefore
$$
\lim_{\alpha\rightarrow +\infty}\limsup_{\l\rightarrow 0\,;\, t\rightarrow +\infty \,;\, \l^2t=\alpha}
\left(\frac{\widehat\E_0^\l [\tau_2 -\tau_1]}{t} \left(\left\lfloor\frac{t}{\widehat\E_0^\l [\tau_2 -\tau_1]}\right\rfloor-1\right) -1\right) =0\,.
$$
We also proved in  (\ref{ellbdd}) that  $\frac{\vert\ell(\l)\vert}{\l}$ remains bounded
for $\l$ tending to $0$. Hence we see that the second term in the right hand side of (\ref{firstisdiff}) also
tends to $0$, thus proving that
(\ref{truelim}) holds true.

To show (\ref{errorterm}), we need the following lemma, whose proof is deferred.

\begin{lemma}\label{lm:tauasy}
We have, for each $\varepsilon > 0$,
\begin{equation}\label{tauVeqt}
\lim_{\alpha\rightarrow +\infty}\limsup_{\l\rightarrow 0\,;\, t\rightarrow +\infty \,;\, \l^2t=\alpha}
\widehat\P_0^\l \left[\vert \tau_{n(t)}-t\vert \geq \varepsilon t \right] =0\, .
\end{equation}
\end{lemma}

We will split the integration in (\ref{errorterm}) according to the partition in the two events
$\{\vert \tau_{n(t)}-t\vert \geq \varepsilon t \}$ and
$\{(1 - \varepsilon)t < \tau_{n(t)} < (1 + \varepsilon)t\}$. We will show that, for each $\varepsilon > 0$,
\begin{equation}\label{split1}
\lim_{\alpha\rightarrow +\infty}\limsup_{\l\rightarrow 0\,;\, t\rightarrow +\infty \,;\, \l^2t=\alpha}
\left\vert \frac{1}{\l t}\widehat\E_0^\l\left[ (X(t) - X(\tau_{n(t)})\1(\vert \tau_{n(t)}-t\vert \geq \varepsilon t )\right]\right\vert =0\, ,
\end{equation}
and that
\begin{equation}\label{split3}
\lim_{\varepsilon \to 0}\limsup_{\alpha\rightarrow +\infty}\limsup_{\l\rightarrow 0\,;\, t\rightarrow +\infty \,;\, \l^2t=\alpha}
\left\vert \frac{1}{\l t}\widehat\E_0^\l\left[ (X(t) - X(\tau_{n(t)})\1((1-\varepsilon)t < \tau_{n(t)}<  (1 + \varepsilon)t)\right]\right\vert =0\, .
\end{equation}

To show (\ref{split1}), we use H\"older's inequality to get that
$$
\left\vert\widehat\E_0^\l\left[ (X(t) - X(\tau_{n(t)})\1(\vert \tau_{n(t)}-t\vert \geq \varepsilon t )\right]\right\vert
$$
$$
\leq \widehat\E_0^\l\left[\left\vert X(t) - X(\tau_{n(t)}\right\vert^2\right]^{1/2}
\widehat\P_0^\l\left[\vert \tau_{n(t)}-t\vert \geq \varepsilon t \right]^{1/2} 
$$
$$
\leq 2 \left(\widehat\E_0^\l\left[\vert X(t)\vert^2\right] + \widehat\E_0^\l\left[\vert X(\tau_{n(t)})\vert^2\right]\right)^{1/2} 
\widehat\P_0^\l\left[\vert \tau_{n(t)}-t\vert \geq \varepsilon t \right]^{1/2}\, ,
$$
(for the last inequality, we used the formula $(x-y)^2 \leq 2(x^2 + y^2)$).
Now, we conclude with (\ref{eq:max****}), (\ref{kdistmoment}) and Lemma \ref{lm:tauasy}.

For (\ref{split3}), note that
$$
\left\vert\widehat\E_0^\l\left[ (X(t) - X(\tau_{n(t)})\1((1-\varepsilon)t < \tau_{n(t)}<  (1 + \varepsilon)t)\right]\right\vert
$$
$$
\leq \E_0^\l\left[\max\limits_{(1-\varepsilon)t \leq u< s \leq (1+\varepsilon)t}\vert X(u) - X(s)\vert\right]
\leq 2  \E_0^\l\left[\max\limits_{(1-\varepsilon)t \leq u \leq (1+\varepsilon)t}\vert X(u) - X((1-\varepsilon)t)\vert\right]\, .
$$
Using the Markov property, we see that the last term equals
$$
g(\l, \varepsilon, t) : = 2  \E_0^\l\left[\max\limits_{u \leq 2\varepsilon t}\vert X(u)\vert\right]\, .
$$
Due to (\ref{eq:max****}), for each $\varepsilon$,
$$
\sup_{\alpha \varepsilon \geq 1} \limsup_{\l\rightarrow 0\,;\, t\rightarrow +\infty \,;\, \l^2\varepsilon t=\alpha}
\frac{g(\l, \varepsilon, t)}{\l \varepsilon t} < \infty\, ,
$$
and this proves (\ref{split3}).

It remains to prove Lemma \ref{lm:tauasy}. We first show
\begin{lemma}\label{lm:LNNtau}
For all $\varepsilon > 0$,
\begin{equation}\label{taubound}
\widehat \P_0^\l \left[\left\vert\tau_k - k\E_0^\l [\tau_2 - \tau_1]\right\vert \geq \varepsilon k \E_0^\l [\tau_2- \tau_1]\right] \to 0\,  \hbox { for } k \to \infty\, ,
\end{equation}
uniformly for $\l \leq 1$.
\end{lemma}

{\it Proof of Lemma \ref{lm:LNNtau}}.\\
We have
\begin{eqnarray}
&  & \widehat\P_0^\l \left[\left\vert\frac{\tau_k}{k} - \widehat\E_0^\l [\tau_2-\tau_1]\right\vert \geq \varepsilon \widehat\E_0^\l [\tau_2 -\tau_1]\right] \cr
&\leq& \frac{1}{\varepsilon^2 \widehat\E_0^\l [\tau_2 -\tau_1]^2} \widehat\E_0^\l \left[\left(\frac{\tau_k}{k} - \widehat\E_0^\l [\tau_2-\tau_1]\right)^2\right]\cr
&\leq& \frac{1}{k^2\varepsilon^2 \widehat\E_0^\l [\tau_2-\tau_1]^2} \widehat\E_0^\l \left[\left(\tau_1 -\widehat\E_0^\l [\tau_2-\tau_1]
+ \sum\limits_{j=2}^k \left(\tau_j -\tau_{j-1} - \widehat\E_0^\l [\tau_2 -\tau_1]\right)\right)^2\right]\cr
&=& \frac{1}{k^2\varepsilon^2 \widehat\E_0^\l [\tau_2-\tau_1]^2}
\left(\widehat\E_0^\l \left[(\tau_1 - \widehat\E_0^\l [\tau_2-\tau_1])^2\right]
+ \sum\limits_{j=2}^k \widehat\E_0^\l \left[\left(\tau_j -\tau_{j-1} - \widehat\E_0^\l [\tau_2 -\tau_1]\right)^2\right]\right)\nonumber
\end{eqnarray}
where we used the independence property stated in Theorem \ref{thm:renew}. But the last term equals
$$
\frac{1}{k^2\varepsilon^2 \widehat\E_0^\l [\tau_2-\tau_1]^2} \left(\widehat\E_0^\l [\tau_1^2] - 2 \widehat\E_0^\l [\tau_1]\widehat\E_0^\l [\tau_2-\tau_1] + \widehat\E_0^\l [(\tau_2-\tau_1)^2]\right)
$$
$$
+ \frac{k-1}{k^2\varepsilon^2 \widehat\E_0^\l [\tau_2-\tau_1]^2} \widehat\E_0^\l \left[\left(\tau_2 -\tau_1 - \widehat\E_0^\l [\tau_2 -\tau_1]\right)^2\right]\, ,
$$
and we conclude, using Lemma \ref{lm:tau}. \qed

Finally, Lemma \ref{lm:tauasy} follows from Lemma \ref{lm:LNNtau} by taking $k= n(t)$: due to (\ref{eq:*}),\\
$ \lim_{\alpha\rightarrow +\infty}\limsup_{\l\rightarrow 0\,;\, t\rightarrow +\infty \,;\, \l^2t=\alpha} n(t) = \infty $ and
\begin{equation}\label{Vasy}
\lim_{\alpha\rightarrow +\infty}\limsup_{\l\rightarrow 0\,;\, t\rightarrow +\infty \,;\, \l^2t=\alpha} \frac{n(t) \widehat\E_0^\l [\tau_2 -\tau_1]}{t} = 1 \, .
\end{equation}
To see that (\ref{Vasy}) holds true, note that
$$
t\left(1 - \frac{\l^2 \widehat\E_0^\l [\tau_2 -\tau_1]}{\l^2 t}\right) \leq
n(t)\widehat\E_0^\l [\tau_2-\tau_1]
\leq t
$$
and use (\ref{eq:*}). \qed


\vskip 1cm
\section{Proof of the Theorem}
\setcounter{equation}{0}
\label{sec:VI}

Combine (\ref{eq:lln}) with Proposition \ref{lm:LR}. \qed


\vskip 1cm
\section{Extension to measurable coefficients}
\setcounter{equation}{0} \label{sec:VII}

Here we explain how the approach developed in the previous Sections
of the paper can be extended to deal with measurable coefficients.
Thus the assumptions in force in this Section are {\bf Assumptions
1, 3 and 4} that remain unchanged and {\bf Assumption 2} is replaced
by the following weaker statement:

\medskip\noindent
{\bf Assumption 2':} for any environment $\o$, the functions
$x\rightarrow V^\o(x)$ and $x\rightarrow \sigma^\o(x)$ are
measurable.

In such generality, it is not possible to use stochastic
differential equations to define the processes $X^\omega$ or
$X^{\l,\,\o}$ anymore so that our first task is to give an
alternative construction.

In the following discussion we fix an environment $\o$ satisfying
{\bf assumptions 1, 2', 3} and {\bf 4}. The case $\l=0$ is included.

Let $p^{\,\l,\,\o}(t,x,y)$ be the kernel associated to the operator
$\LL^{\l,\,\o}$ in equation (\ref{generator}), or, in the case
$\l=0$ to the operator $\LL^\o$ in equation (\ref{eq:gen}) now
understood in weak distributional sense and let $T^{\l,\,\o}_t$ be
the corresponding semigroup
$$T^{\l,\,\o}_t f(x)=\int p^{\,\l,\,\o}(t,x,y) f(y) dy\,.$$
It follows from Aronson's estimate (see \cite{kn:aronson}) and
Harnack's inequality (see \cite{kn:GT}) that $T^{\l,\,\o}_t$ maps
continuous functions  vanishing  at infinity to continuous functions
vanishing  at infinity. For such functions we also have
$\lim_{t\rightarrow 0} T^{\l,\,\o}_t f(x)=f(x)$ for all $x$. Thus
$T^{\l,\,\o}_t$ is a Feller semigroup and it follows from
\cite{kn:EK} Theorem 2.7 that there exists a conservative Hunt
process with continuous paths whose semigroup is $T^{\l,\,\o}$. We
denote its law on path space $C(\R_+,\R^d)$ with $P_x^{\l,\,\o}$ and
$E_x^{\l,\,\o}$ for the corresponding expectation. Observe that
$P_x^{\l,\,\o}[X(0)=x]=1$ for all $x$.

It is proved in \cite{kn:lejay} Proposition 1 that, for almost all
$\o$'s, under $P_0^{\o}$, the canonical process satisfies an
invariance principle with some effective diffusivity matrix
$\Sigma$. It follows from Aronson's estimate that $\Sigma$ is also
the asymptotic covariance so that Definition \ref{df:Sigma} goes
through.

On the other hand all statements in Section \ref{sec:IV} as well as
Lemma \ref{lm:btrans} in Section \ref{sec:V} were proved under {\bf
Assumption 2'} only. Thus the construction of regeneration times we
gave in Section \ref{sec:V} and all the upper bounds on $\tau_1$ are
still valid with measurable coefficients. Therefore we see that the
law of large numbers is satisfied $\Q$ almost surely under
$P_0^{\l,\,\o}$ for all positive $\l$ and Definition \ref{df:velo}
can still be used as the definition of the effective drift.

Having defined the effective diffusivity matrix and the effective
drift we claim that the Einstein relation stated in Theorem
\ref{theo} holds true with {\bf Assumption 2} being replaced by {\bf
Assumption 2'}.

The main difficulties in extending the proofs of the previous
sections to measurable coefficients appear in justifying the
Girsanov transform and time change arguments from Section
\ref{sec:III}. Following \cite{kn:lejay}, in order to do it we shall
appeal to Dirichlet form theory, as exposed in \cite{kn:FOT}, and
related stochastic calculus for Dirichlet processes. Observe that a
direct application of Dirichlet form theory a priori only provides
information under $P_x^{\l,\,\o}$ for all $\o$ but only for
(Lebesgue) almost all $x$ (in fact for quasi all $x$ but we won't
use fine topological notions here) and therefore, as a consequence
of the translation invariance of $\Q$, under $P_0^{\l,\,\o}$ for
almost all $\o$. Therefore most claims in Sections \ref{sec:III},
such as Lemma \ref{lm:ap1} or formula (\ref{eq:gir}), should now be
understood 'for $\Q$ almost all $\o$'s'. We let the reader convince
herself that this does not affect the proofs.

We use the notation $\rho_\l(x)=e^{{\hat \l}\cdot x}$ for $\l\geq
0$. Let $L^2(\rho_\l^2)$ be the space of square integrable functions
with respect to the measure $\rho_\l^2(x)\,dx$. Define
$H_1(\rho_\l^2)$ to be the space of functions in $L^2(\rho_\l^2)$
whose gradient is also square integrable with respect to the measure
$\rho_\l^2(x)\,dx$. Let \beqn\label{eq:exprdir}{\cal
E}^{\l,\,\o}(f,f):=\frac 12 \int \vert\sigma^\o(x)\nabla
f(x)\vert^2\, e^{-2V^\o(x)}\,\rho_\l^2(x)\,dx\,.\eeqn Then $({\cal
E}^{\l,\,\o},H_1(\rho_\l^2))$ is a regular Dirichlet form. We claim
that

\begin{lemma}
$({\cal E}^{\l,\,\o},H_1(\rho_\l^2))$ is the Dirichlet form of the
semigroup $T^{\l,\,\o}$ on $L^2(e^{2{\hat \l}\cdot x-2V^\o(x)}dx)$.
\end{lemma}

(Note that this fact is already used in \cite{kn:lejay} but without
justification.)

{\it Proof}

We first observe that $T^{\l,\,\o}$ is indeed a strongly continuous
symmetric semigroup on $L^2(\rho_\l^2)$.

Let $t>0$ and define the approximating bilinear forms
\beqn\label{eq:app} {\cal E}^{t,\,\l,\,\o}(f,f):=\frac 1t \int
(f(x)-T^{\l,\,\o}_t f(x))f(x)\,
e^{-2V^\o(x)}\,\rho_\l^2(x)\,dx\,.\eeqn  A function $f$ belongs to
the domain   of the Dirichlet form associated to the semigroup
$T^{\l,\,\o}$ if and only if ${\cal E}^{t,\,\l,\,\o}(f,f)$ is
bounded in $t$ and the limit as $t$ tends to $0$ is then the value
of the form, see Lemma 1.3.4 in \cite{kn:FOT}.

A straighforward integration by parts shows that, as $t$ tends to
$0$, then ${\cal E}^{t,\,\l,\,\o}(f,f)$ converges to ${\cal
E}^{\l,\,\o}(f,f)$ on the $L^2$ domain of the generator
$\LL^{\l,\,\o}$, say ${\cal D}^{\l,\,\o}$. Since the function
$t\rightarrow {\cal E}^{t,\,\l,\,\o}(f,f)$ is decreasing, it implies
that $${\cal E}^{t,\,\l,\,\o}(f,f)\leq {\cal E}^{\l,\,\o}(f,f)$$ for
$f\in {\cal D}^{\l,\,\o}$. This inequality extends by density to all
functions in $H_1(\rho_\l^2)$. Thus we have proved that the
Dirichlet form of the semigroup $T^{\l,\,\o}$ is well defined and
coincides with ${\cal E}^{\l,\,\o}$ on $H_1(\rho_\l^2)$ or, in other
words, that it is an extension of the form $({\cal
E}^{\l,\,\o},H_1(\rho_\l^2))$. But since $({\cal
E}^{\l,\,\o},H_1(\rho_\l^2))$ is its own maximal Markovian
extension, see Theorem 3.3.1 in \cite{kn:FOT}, both forms coincide.
\qed

>From now on we will drop the superscript $\l$ from the notation when
$\l=0$. We now consider properties of the canonical process $X$ for
a fixed environment $\o$ and under $P_x^{\o}$ for almost every
starting point $x$. The function $\phi(x)=e_1\cdot x$ locally
belongs to the domain of the Dirichlet form ${\cal E}^{\o}$. From
Theorem 5.5.1 in \cite{kn:FOT} we deduce that $e_1\cdot
X(t)-e_1\cdot X(0)$ is a local Dirichlet process under $P_x^{\o}$
for almost every starting point $x$. Thus $e_1\cdot X(t)-e_1\cdot
X(0)$ admits a unique Fukushima decomposition as the sum of a local
martingale, say $\bar B$, and a process of locally vanishing
quadratic variation. The bracket of $\bar B$ is given by Theorem
5.5.2 in \cite{kn:FOT} and satisfies \beqn\label{eq:bracket}
\langle\bB\rangle(t)=\int_0^t \vert \sigma^\o(X(s))\, e_1 \vert^2 \,
ds \,.\eeqn By {\bf assumption 3}, we have $\langle\bB\rangle(t)\leq
\kappa^{-1}t$. Thus we see that in fact $\bB$ is a square integrable
martingale. By the same argument the exponential local martingale
$e^{\mu\bB (t) -\frac {\mu^2}2 \langle\bB\rangle(t)}$ is also seen
to be a martingale for all $\mu$ and one proves as in Lemma
\ref{lm:ap1} that \beqn\label{eq:expa} E\left[\, e^{\mu\l\bB (t)
-\frac {\mu\l^2}2 \langle\bB\rangle(t)}\right] \leq e^{(\mu-1)\frac
{\mu\l^2}2\frac t{\kappa}}\,, \eeqn for every $\o$ and almost all
$x$. Note that the translation invariance of $\Q$ then implies that
(\ref{eq:expa}) also holds for almost all $\o$ with $x=0$.

We now justify the Girsanov formula:

\begin{proposition}\label{prop:asgir}
For any environment $\o$, for almost any $x$, any $t$ and any
continuous bounded function $F$ we have
$$E^{\l,\,\o}_x[F(X([0,t]))]=E^\o_x[F(X([0,t]))\,
e^{\l\bB (t) -\frac {\l^2}2 \langle\bB\rangle(t)}]\,.$$
\end{proposition}

The translation invariance of $\Q$ then implies that for almost all
environments
$$E^{\l,\,\o}_0[F(X([0,t]))]=E^\o_0[F(X([0,t]))\, e^{\l\bB (t)
-\frac {\l^2}2 \langle\bB\rangle(t)}]\,,$$ for any continuous and
bounded function $F$ and we have obtained the almost sure version of
the Girsanov formula (\ref{eq:gir}) that is sufficient to proceed
through the proofs of the previous Sections.

{\it Proof}

We would like to invoke Theorem 3.1 of \cite{kn:Fitz} but
unfortunately $\phi$ does not belong to $L^2$ so that some work is
needed.

Let $G_n$ be the ball centered at the origin with radius $n$, and
let $\zeta_n$ be the exit time from $G_n$.

We shall first prove that

\beqn \label{eq:locgir} E^{\l,\,\o}_x[F(X([0,t]))\,;\,
t<\zeta_n]=E^\o_x[F(X([0,t]))\, e^{\l\bB (t) -\frac {\l^2}2
\langle\bB\rangle(t)}\,;\, t<\zeta_n]\,.\eeqn

Using the bound (\ref{eq:expa}), it is then possible to let $n$ tend
to infinity and deduce Proposition \ref{prop:asgir} from
(\ref{eq:locgir}).

Choose a function $\phi_n$ that coincides with $\phi$ on $G_n$, is
smooth and has compact support. Let $M_n$ be the martingale part of
the process $\phi_n(X(t))-\phi_n(X(0))$ in its Fukushima
decomposition under $P_x^{\o}$, and let $Z_n(t):=e^{\l M_n(t) -\frac
{\l^2}2 \langle M_n\rangle(t)}$.

Define
$$ Q^{\,n,\,\l,\,\o}_t f(x):=E^\o_x[f(X(t))Z_n(t)]\,.$$
Then $Q^{\,n,\,\l,\,\o}_t$ defines a strongly continuous Markovian
semigroup on $L^2(dx)$. We need the following

\begin{lemma}
The Dirichlet form of the semigroup $Q^{n\,,\,\l,\,\o}_t$ acting on
$L^2(e^{2\phi_n(x)-2V^\o(x)}dx)$ is
$$\frac 12 \int \vert\sigma^\o(x)\nabla
f(x)\vert^2\, e^{2\phi_n(x)-2V^\o(x)}\,dx\,,$$ with domain $H_1(dx)$.
\end{lemma}

{\it Proof}

The Revuz measure of the positive continuous additive functional
$\langle M_n \rangle$ is $\vert\nabla \sigma^\o(x)\phi_n(x)\vert^2
dx$ which is easily seen to belong to the Hardy class since the
gradient of $\phi_n$ is uniformly bounded. Besides condition (3.8)
in \cite{kn:Fitz} is fulfilled whenever $\nabla\phi_n$ is uniformly
bounded. Thus Theorem 3.1 of \cite{kn:Fitz} applies. We get that the
quadratic form of the semigroup $Q^{n\,,\,\l,\,\o}_t$ acting on
$L^2(e^{-2V^\o(x)}dx)$ is \beqn  {\cal Q}_n(f,g):= &&\frac 12\int
\sigma^\o(x)\nabla f(x)\cdot\sigma^\o(x)\nabla g(x)\,
e^{-2V^\o(x)}\,dx\cr \label{eq:q}&-&\int g(x)\,\sigma^\o(x)\nabla
f(x)\cdot\sigma^\o(x)\nabla \phi_n(x)\, e^{-2V^\o(x)}\,dx\,, \eeqn
see formula (3.3) in \cite{kn:Fitz}.

We now use the same approximating sequence as in (\ref{eq:app}); see
also the explanation on page 242 of \cite{kn:Fitz}. We then  know
that for all functions $f,g\in H_1(dx)$
$${\cal Q}_n(f,g)=\lim_{t\rightarrow0}
\frac 1t \int (f(x)-Q^{n\,,\,\l,\,\o}_t f(x))g(x)\,
e^{-2V^\o(x)}\,dx\,.$$ Applying this formula to the function
$ge^{2\phi_n}$ and using (\ref{eq:q}), we deduce that the Dirichlet
form of the semigroup $Q^{n\,,\,\l,\,\o}_t$ acting on
$L^2(e^{2\phi_n(x)-2V^\o(x)}dx)$ is \beqnn&&\lim_{t\rightarrow0}
\frac 1t \int (f(x)-Q^{n\,,\,\l,\,\o}_t f(x))g(x)e^{2\phi_n(x)}\,
e^{-2V^\o(x)}\,dx\\
&=&\frac 12\int \sigma^\o(x)\nabla f(x)\cdot\sigma^\o(x)\nabla
(g(x)e^{2\phi_n(x)})\, e^{-2V^\o(x)}\,dx\\ &-&\int
g(x)\,\sigma^\o(x)\nabla f(x)\cdot\sigma^\o(x)\nabla \phi_n(x)\,
e^{2\phi_n(x)-2V^\o(x)}\,dx \\
&=& \frac 12 \int \sigma^\o(x)\nabla f(x)\cdot\sigma^\o(x)\nabla
g(x)\, e^{2\phi_n(x)-2V^\o(x)}\,dx\,.\eeqnn \qed

We can now conclude the proof of Proposition \ref{prop:asgir}.

Comparing the expression of the Dirichlet form generated by
$Q^{n\,,\,\l,\,\o}_t$ we just obtained with formula
(\ref{eq:exprdir}) for ${\cal E}^{\l,\,\o}$ and observing that
$\phi_n=\phi$ on $G_n$ we see that the parts of both these Dirichlet
forms on $G_n$ coincide. Clearly $\bar B$ and $M_n$ also coincide up
to time $\zeta_n$. Thus we obtain (\ref{eq:locgir}).\qed

Finally we should say a word about the time change argument used in
the proof of Proposition \ref{lm:LR} in the case $V\not=0$. Define
$Y^{\l,\,\o}$ to be the Hunt process with Dirichlet form $({\cal
E}^{\l,\,\o},H_1(\rho_\l^2))$ with reference measure
$\rho_\l^2(x)dx$. Theorem 6.2.1 in \cite{kn:FOT} implies that the
process obtained by time changing $Y^{\l,\,\o}$ through the additive
functional $ A^{\l,\,\o}(t)=\int_0^t e^{-2V^\o(Y^{\l,\,\o}(s))}\,
ds$ admits as Dirichlet form $({\cal E}^{\l,\,\o},H_1(\rho_\l^2))$
with reference measure $e^{-2V^\o(x)}\rho_\l^2(x)dx=e^{2{\hat
\l}\cdot x-2V^\o(x)}$ so that, for almost any initial point $x$, the
law of $Y^{\l,\,\o}\circ(A^{\l,\,\o})^{-1}$ coincides with
$P^{\l,\,\o}_x$.





\vskip 1cm
{\it Acknowledgement: we thank J.C. Mourrat for his careful reading of a preliminary version of the paper.}


\frenchspacing
\bibliographystyle{plain}

\end{document}